\documentclass[12pt]{amsart}

\calclayout

\usepackage{amssymb}
\usepackage{latexsym}
\usepackage{exscale}
\usepackage{graphicx}

\usepackage[colorlinks=true, pdfstartview=FitV, linkcolor=red,
citecolor=blue, urlcolor=blue]{hyperref}

\usepackage{color}

\allowdisplaybreaks

\numberwithin{equation}{section}

\theoremstyle{plain}
\newtheorem{theorem}[equation]{Theorem}
\newtheorem{lemma}[equation]{Lemma}
\newtheorem{cor}[equation]{Corollary}

\theoremstyle{definition}

\theoremstyle{remark}
\newtheorem{remark}[equation]{Remark}

\newcommand{\re}{\mathbb{R}}
\newcommand{\rn}{{\mathbb{R}^n}}

\newcommand{\F}{\mathcal{F}}
\newcommand{\R}{\mathcal{R}}

\DeclareMathOperator*{\esssup}{ess\,sup}

\def\Xint#1{\mathchoice
  {\XXint\displaystyle\textstyle{#1}}%
  {\XXint\textstyle\scriptstyle{#1}}%
  {\XXint\scriptstyle\scriptscriptstyle{#1}}%
  {\XXint\scriptscriptstyle\scriptscriptstyle{#1}}%
  \!\int}
\def\XXint#1#2#3{{\setbox0=\hbox{$#1{#2#3}{\int}$}
    \vcenter{\hbox{$#2#3$}}\kern-.5\wd0}}

\def\avgint{\Xint-}

\title[Limited range multilinear extrapolation]{Limited range
  multilinear extrapolation with applications to the bilinear Hilbert
  transform}

\author{David Cruz-Uribe, OFS}
\address{David Cruz-Uribe, OFS \\ 
Department of Mathematics \\ 
University of Alabama \\
Tus\-ca\-loosa, AL 35487, USA}
\email{dcruzuribe@ua.edu}

\author{Jos\'e Mar{\'\i}a Martell}
\address{Jos\'e Mar{\'\i}a Martell
\\
Instituto de Ciencias Matem\'aticas CSIC-UAM-UC3M-UCM
\\
Consejo Superior de Investigaciones Cient{\'\i}ficas
\\
C/ Nicol\'as Cabrera, 13-15
\\
E-28049 Madrid, Spain} \email{chema.martell@icmat.es}

\thanks{The first author is supported by NSF Grant DMS-1362425 and
  research funds from the Dean of the College of Arts \& Sciences,
  University of Alabama.  The second author acknowledges financial
  support from the Spanish Ministry of Economy and Competitiveness,
  through the ``Severo Ochoa'' Programme for Centres of Excellence in
  R\&D” (SEV-2015-0554). He also acknowledges that the research
  leading to these results has received funding from the European
  Research Council under the European Union's Seventh Framework
  Programme (FP7/2007-2013)/ ERC agreement no. 615112 HAPDEGMT. The
  authors would like to thank Francesco di Plinio for suggesting this
  problem to us and for helpful discussions about the bilinear Hilbert
  transform.  The authors would also like to thank Sheldy Ombrosi for
  suggesting the application to Marcinkiewicz-Zygmund estimates. The authors express their gratitude to
  Camil Muscalu and Cristina Benea  for helpful discussions about the vector-valued inequalities for the bilinear Hilbert transform.	
	Finally the first author would  like to thank the second
  author for his hospitality during two visits to Madrid where much of
  the work on the project was done.}



\date{\today}

\begin{document}

\begin{abstract}
We prove a limited range, off-diagonal extrapolation theorem that
generalizes a number of results in the theory of Rubio de Francia
extrapolation, and use this to prove a limited range, multilinear
extrapolation theorem.   We give two applications of this result to
the bilinear Hilbert transform.  First, we give sufficient conditions on a
pair of weights $w_1,\,w_2$ for the bilinear Hilbert transform to
satisfy weighted norm inequalities of the form
\[ BH : L^{p_1}(w_1^{p_1}) \times L^{p_2}(w_2^{p_2}) 
\longrightarrow L^p(w^p), \]
where $w=w_1w_2$ and
$\frac{1}{p}=\frac{1}{p_1}+\frac{1}{p_2}<\frac{3}{2}$.  This improves
the recent results of Culiuc {\em et al.}~by increasing the families
of weights for which this inequality holds and by pushing the lower
bound on $p$ from $1$ down to $\frac{2}{3}$, the critical index from
the unweighted theory of the bilinear Hilbert transform.  Second,
as an easy consequence of our method we obtain that the bilinear
Hilbert transform satisfies some vector-valued inequalities with
Muckenhoupt weights. This reproves and generalizes some of the
vector-valued estimates obtained by Benea and Muscalu in the
unweighted case. We also generalize recent results of Carando,
{\em et al.} on Marcinkiewicz-Zygmund estimates for multilinear
Calder\'on-Zygmund operators.
\end{abstract}
\maketitle

\section{Introduction}
\label{section:intro}

The Rubio de Francia theory of extrapolation is a powerful 
tool in harmonic analysis.  In its most basic form, it shows that if,
for a fixed value $p_0$, $1<p_0<\infty$, an operator $T$ satisfies a
weighted norm inequality of the form
\begin{equation} \label{eqn:hyp-RdF}
 \|Tf\|_{L^{p_0}(w)} \leq C\|f\|_{L^{p_0}(w)} 
\end{equation}
for every weight $w$ in the Muckenhoupt class $A_{p_0}$, then for every
$p$, $1<p<\infty$, 
\begin{equation} \label{eqn:con-RdF}
  \|Tf\|_{L^{p}(w)} \leq C\|f\|_{L^{p}(w)} 
\end{equation}
whenever $w\in A_p$.  Since its discovery in the early 1980s,
extrapolation has been generalized in a variety of ways, yielding
weak-type inequalities, vector-valued inequalities, and inequalities
in other scales of Banach function spaces.  We refer the reader
to~\cite{cruz-martell-perezBook} for the development of extrapolation;
for more recent results we refer the reader to~\cite{
  Cruz-Uribe2016Extrapolation-a, CruzUribe:2015km, MR2754896}.

Extrapolation has been also extended to the multilinear
setting.  In~\cite{grafakos-martell04} it was shown that if a given operator $T$ satisfies
\[ \|T(f_1,\ldots,f_m)\|_{L^p((w_1\cdots w_m)^p)}
\leq C\prod_{j=1}^m \|f_j\|_{L^{p_j}(w_j^{p_j})} \]
for fixed exponents $1<p_1,\ldots,p_m<\infty$, 
$\displaystyle \frac{1}{p} = \sum_{j=1}^m \frac{1}{p_j}$,
and all weights $w_j^p \in A_{p_j}$, then the same estimate holds
for all possible values of $p_j$.   An extension to the scale of
variable Lebesgue spaces was given in~\cite{CruzUribe:2016wv}. 

\medskip 

In this paper we develop a theory of limited range, multilinear
extrapolation.  In the linear case, limited range extrapolation was
developed in~\cite{auscher-martell07} by Auscher and the second
author.  They proved that if
inequality~\eqref{eqn:hyp-RdF} holds for a given $0<p_-<p_0<p_+<\infty$ and for all
$w\in A_{\frac{p_0}{p_-}}\cap RH_{\big(\frac{p_+}{p_0}\big)'}$, then for all
$p_-<p<p_+$ and $w\in A_{\frac{p}{p_-}}\cap RH_{\big(\frac{p_+}{p}\big)'}$,
\eqref{eqn:con-RdF} holds.  Conditions like this arise naturally in
the study of the Riesz transforms and other operators associated to
elliptic differential operators.

Our first theorem extends limited range extrapolation to the
multilinear setting.    To state our results we use the abstract
formalism of extrapolation families.    Given $m\geq 1$, hereafter $\F$ will
denote a family of $(m+1)$-tuples $(f,f_1,\ldots,f_m)$ of non-negative
measurable functions.   This approach to extrapolation has the
advantage that, for instance, vector-valued inequalities are an
immediate consequence of our extrapolation results.  We will discuss
applying this formalism to prove norm inequalities for specific operators
below.  For complete discussion of this
approach to extrapolation in the linear setting, see~\cite{cruz-martell-perezBook}.

\begin{theorem} \label{thm:extrapol-multi}
Given $m\ge 1$, let $\F$ be a family of
extrapolation $(m+1)$-tuples.  
For each $j$, $1\leq j\leq m $, suppose we have parameters $r_j^-$ and $r_j^+$,
and an exponent $p_j\in (0,\infty)$, $0\le  r_j^-\le p_j \le r_j^+\le \infty$, such that given any collection of weights
$w_1,\ldots,w_m$ with $w_j^{p_j} \in A_{\frac{p_j}{r_j^-}}\cap
RH_{\big(\frac{r_j^+}{p_j}\big)'}$ and $w=w_1\cdots w_m$, we have the inequality
\begin{equation}
\label{eq:extrap-multi-hyp}
\|f\|_{L^p(w^p)} 
\leq
C\prod_{j=1}^m \|f_j\|_{L^{p_j}(w_j^{p_j})}
\end{equation}
for all $(f,f_1,\ldots,f_m)\in \F$  such that $\|f\|_{L^p(w^p)}<\infty$,  where $\displaystyle \frac{1}{p} = \sum_{j=1}^m \frac{1}{p_j}$
and $C$ depends on
$n,\,p_j,\,[w_j]_{A_{\frac{p_j}{r_j^-}}},\,[w_j]_{RH_{\big(\frac{r_j^+}{p_j}\big)'}}$. 
Then for all exponents $q_j$, $r_j^-<q_j<r_j^+$,  all weights $w_j^{q_j} \in
A_{\frac{q_j}{r_j^-}}\cap RH_{\big(\frac{r_j^+}{q_j}\big)'}$ and $w=w_1\cdots w_m$,
\begin{equation}\label{eq:extrap-multi-con}
\|f\|_{L^q(w^q)} 
\leq
C\prod_{j=1}^m \|f_j\|_{L^{q_j}(w_j^{q_j})}, 
\end{equation}
for all $(f,f_1,\ldots,f_m)\in \F$  such that $\|f\|_{L^q(w^q)}<\infty$, where 
$\displaystyle \frac{1}{q} = \sum_{j=1}^m \frac{1}{q_j}$
and $C$ depends on
$n,\,p_j,\,q_j,\,[w_j]_{A_{\frac{q_j}{r_j^-}}},\,[w_j]_{RH_{\big(\frac{r_j^+}{q_j}\big)'}}$. 
Moreover, for the same family of exponents and weights, and for all exponents $s_j$, $r_j^-<s_j<r_j^+$,
\begin{equation}\label{eq:extrap-multi-con-vv}
\left\|\bigg(\sum_k (f^k)^s\bigg)^{\frac1s}\right\|_{L^q(w^q)} 
\leq
C\prod_{j=1}^m 
\left\|\bigg(\sum_k (f_j^k)^{s_j}\bigg)^{\frac1{s_j}}\right\|_{L^{q_j}(w_j^{q_j})}, 
\end{equation}
for all $\big\{(f^k,f_1^k,\ldots,f_m^k)\}_k\subset \F$  such that the left-hand side is finite and where  
$\displaystyle \frac{1}{s} = \sum_{j=1}^m \frac{1}{s_j}$
and $C$ depends on
$n,\,p_j,\,q_j,\,s_j,\,[w_j]_{A_{\frac{q_j}{r_j^-}}},\,[w_j]_{RH_{\big(\frac{r_j^+}{q_j}\big)'}}$. 

\end{theorem}

\begin{remark}
When $r_j^-=1$ and $r_j^+=\infty$ in Theorem~\ref{thm:extrapol-multi}
we get a version of the multilinear
extrapolation theorem from~\cite{grafakos-martell04} for extrapolation
families.  The original result was given in terms of operators.
\end{remark}

\medskip

Theorem~\ref{thm:extrapol-multi} is a consequence of a linear, restricted range, off-diagonal
extrapolation theorem, which we believe is of interest in its own
right.  It generalizes the classical Rubio de Francia extrapolation,
the off-diagonal extrapolation theory of Harboure, Mac{\'\i}as and
Segovia~\cite{Harboure:1988vi}, and the limited range extrapolation
theorem proved by Auscher and the second
author~\cite{auscher-martell07}.

\begin{theorem} \label{thm:off-diag-limited}
Given $0\le p_-<p_+\leq \infty$ and a family of extrapolation pairs $\F$, suppose that for some $p_0,\,q_0\in (0,\infty)$ such that
$p_-\le p_0\le p_+$,  $\frac1{q_0}-\frac1{p_0}+\frac1{p_+}\ge 0$, and all weights $w$ such that $w^{p_0}\in A_{\frac{p_0}{p_-}}
\cap RH_{\big(\frac{p_+}{p_0}\big)'}$, 
\begin{equation} \label{eqn:off-diag-limited1}
\left( \int_\rn f^{q_0}w^{q_0}\,dx \right)^{\frac{1}{q_0}}
\leq C\left( \int_\rn g^{p_0} w^{p_0}\,dx\right)^{\frac{1}{p_0}}
\end{equation}
for all $(f,g) \in \F$ such that $\|f\|_{L^{q_0}(w^{q_0})}<\infty$,
and the constant $C$ depends on
$n,\,p_0,\,q_0,\,[w^{p_0}]_{A_{\frac{p_0}{p_-}}},\,
[w^{p_0}]_{RH_{\big(\frac{p_+}{p_0}\big)'}}$.  Then for every $p$, $q$
such that $p_-<p<p_+$, $0<q<\infty$ and
$\frac{1}{p}-\frac{1}{q}=\frac{1}{p_0}-\frac{1}{q_0}$, and every
weight $w$ such that
$w^{p}\in A_{\frac{p}{p_-}} \cap RH_{\big(\frac{p_+}{p}\big)'}$,
\begin{equation} \label{eqn:off-diag-limited2}
\left( \int_\rn f^{q}w^{q}\,dx \right)^{\frac{1}{q}}
\leq C\left( \int_\rn g^{p} w^{p}\,dx\right)^{\frac{1}{p}}
\end{equation}
for all 
$(f,g) \in \F$ such that $\|f\|_{L^{q}(w^{q})}<\infty$, and $C$
depends on  $n$, $p$, $q$, $[w^p]_{A_{\frac{p}{p_-}}}$,
$[w^p]_{RH_{\big(\frac{p_+}{p}\big)'}}$.  
\end{theorem}

In Theorems~\ref{thm:extrapol-multi}
and~\ref{thm:off-diag-limited} we make the {\em a priori} assumption
that the left-hand sides of both our hypothesis and conclusion are
finite, and this plays a role in the proof.   In certain applications
this assumption is reasonable:  for instance, when proving
Coifman-Fefferman type inequalities
(cf.~\cite{cruz-martell-perezBook}).  
However, when using extrapolation to prove norm
inequalities for operators we would like to remove this assumption, as
the point is to conclude that the left-hand side is finite.  But in
fact, we can do this by an easy approximation argument.  This
immediately yields the following corollaries.

\begin{cor} \label{cor:extrapol-mult} 
Under the same hypotheses as
  Theorem~\ref{thm:extrapol-multi}, if we assume
  that~\eqref{eq:extrap-multi-hyp} holds for all
  $(f,f_1,\ldots,f_m)\in\F$ (whether or not the left-hand side is
  finite) then the conclusion~\eqref{eq:extrap-multi-con} holds for
  all $(f,f_1,\ldots,f_m)\in\F$ (whether or not the left-hand side is
  finite).  Analogously, the vector-valued inequality~\eqref{eq:extrap-multi-con-vv} holds for all
	families $\big\{(f^k,f_1^k,\ldots,f_m^k)\big\}_k\subset\F$ (whether or not the left-hand side is
  finite).
\end{cor}

\begin{cor}\label{corol:extrapol-rest}
  Under the same hypotheses as Theorem \ref{thm:off-diag-limited}, if
 we assume that~\eqref{eqn:off-diag-limited1} holds for all
  $(f,g)\in\F$ (whether or not the left-hand side is finite) then the
  conclusion \eqref{eqn:off-diag-limited2} holds for all $(f,g)\in\F$
  (whether or not the left-hand side is finite).
\end{cor}

\medskip

In the statement of Theorem~\ref{thm:off-diag-limited} there are some
restrictions on the allowable exponents $p$ and $q$.  We make these
explicit here; these restrictions will play a role in the proof
below.  

\begin{remark}\label{remark:rest-expon}
 Define $q_\pm$ by
 \begin{equation}   \label{eq:defqpm}
 \frac{1}{q_\pm}-\frac{1}{p_\pm} = \frac{1}{q_0}-\frac{1}{p_0}.
 \end{equation}
 Because of our assumptions that
 $\frac1{q_0}-\frac1{p_0}+\frac1{p_+}\ge 0$ and $0\le p_-\le p_0\le p_+\le \infty$ it follows that
 $0\le q_-\le q_0\le q_+\leq \infty$. Moreover, the fact that $p_-<p<p_+$ yields that $q_-<q<q_+$.     Note that if we were to allow that $\frac1{q_0}-\frac1{p_0}+\frac1{p_+}<0$, we could choose $p$ very close to $p_+$ and the associated $q$ would be negative, which would not
 make sense.   

Moreover, we have that the following hold:
\begin{enumerate}
\renewcommand{\theenumi}{\roman{enumi}}
\item  If $q_0=p_0$, then $q_{\pm}=p_\pm$ and $q=p$. 

\item If $p_0>q_0$, then
 $0\leq q_-<p_-$, $q_+<p_+\le \infty$ and $q<p$.

\item  If $p_0<q_0$, then
 $0\le p_-<q_-$, $p_+<q_+\le\infty$ and $p<q$. 
\end{enumerate}
\end{remark}

\begin{remark}\label{remark:rest-expon:0}
  When $p_0\ge q_0$ we automatically have that
  $\frac1{q_0}-\frac1{p_0}+\frac1{p_+}\ge 0$.  Further, this implies
  that all of the weights which appear in both our hypothesis and
  conclusion (i.e, $w^{p_0}$, $w^{q_0}$, $w^p$, $w^q$) are in
  $A_\infty$.    Consequently, they are locally integrable, and  so all the
  Lebesgue spaces that appear in the statement contain the
  characteristic functions of compact sets.  In fact, since
  $w^{p_0}\in A_\infty$,  $w^{q_0}\in A_\infty$ (see
  Lemma~\ref{lemma:cjn} below).  The same is
  true for $w^p$ and $w^q$, since by Remark~\ref{remark:rest-expon},
  $p\geq q$.

  When $p_0<q_0$, the condition
  $\frac1{q_0}-\frac1{p_0}+\frac1{p_+}\ge 0$ imposes an upper bound
  for $q_0$: $q_0\le p_0(p_+/p_0)'$.  A similar bound holds for $q$.
  Thus (by Lemma~\ref{lemma:cjn}) $w^{q_0},\, w^q\in A_\infty$ and
  so again all the weights involved are  in $A_\infty$ and thus locally
  integrable. 
\end{remark}

\medskip

Theorem~\ref{thm:off-diag-limited} and Corollary~\ref{corol:extrapol-rest} generalize several known extrapolation
results.
\begin{list}{$(\theenumi)$}{\usecounter{enumi}\leftmargin=.8cm
\labelwidth=.8cm\itemsep=0.2cm\topsep=.1cm
\renewcommand{\theenumi}{\roman{enumi}}}

\item The classical Rubio de Francia extrapolation theorem (see
  e.g.~\cite[Theorems 1.4 and 3.9]{cruz-martell-perezBook} for the
  precise formulation) corresponds to the case $p_-=1$, $p_+=\infty$,
  $q_0=p_0$. 

\item The $A_\infty$ extrapolation theorem
  in~\cite{cruz-uribe-martell-perez04} (see
  also~\cite[Corollary~3.15]{cruz-martell-perezBook}) corresponds to the
  case $p_-=0$, $p_+=\infty$, and $q_0=p_0$.  

\item The extrapolation theorem for weights in the reverse H\"older classes
  \cite[Lemma 3.3, (b)]{martell-prisuelos} corresponds to the case
$p_-=0$, $p_+=1$, and $q_0=p_0$. 

\item The limited range extrapolation theorem in~\cite[Theorem
  4.9]{auscher-martell07} (see also  \cite[Theorems 3.31]{cruz-martell-perezBook}),  corresponds to the case $0<p_-<p_+\le \infty$, $q_0=p_0$. 

\item The off-diagonal extrapolation theorem in \cite{Harboure:1988vi}
  (see  also~\cite[Theorem~3.23]{cruz-martell-perezBook})  corresponds to
  the case $p_-=1$, $p_0<q_0$,
  $p_+=\big(\frac1{p_0}-\frac1{q_0}\big)^{-1}$.  To see this, we
  recall the well-known fact that $w\in A_{p_0, q_0}$, that is,
\[ \sup_Q  \left( \avgint_Q w^{q_0}\,dx\right)^{\frac{1}{q_0}}
\left( \avgint_Q w^{-p_0'}\,dx\right)^{\frac{1}{p_0'}}<\infty, \]
if and only if
$w^{p_0}\in A_{p_0}\cap RH_{\frac{q_0}{p_0}}= A_{\frac{p_0}{p_-}} \cap
RH_{\big(\frac{p_+}{p_0}\big)'}$.  Note that in this case
$\frac1{q_0}-\frac1{p_0}+\frac1{p_+}=0$.
\end{list}

\medskip

Our generalization of  off-diagonal extrapolation involves weighted
norm inequalities that have already appeared in the literature in the context of
  fractional powers of second divergence form elliptic operators with
  complex bounded measurable coefficients.   More precisely, in
  \cite{auscher-martell08} it was shown that for a certain operator
  $T_\alpha$, there exist 
  $1\le r_-<2<r_+\le\infty$ such that  $ T_\alpha : L^r(w^r)
  \rightarrow L^s(w^s)$ for every $r_-<r<s<r_+$ and for
  every
  $w\in A_{1+\frac1{r_-}-\frac1r}\cap
  RH_{s\big(\frac{r_+}{s}\big)'}$.  
By applying Theorem~\ref{thm:off-diag-limited} we could prove
  the same result via extrapolation if we could show that there exists
  $r_-<r_0<s_0<r_+$ such that $T_\alpha : L^{r_0}(w^{r_0}) \rightarrow
  L^{s_0}(w^{s_0})$ for every
  $w\in A_{1+\frac1{r_-}-\frac1{r_0}}\cap
  RH_{s_0(\frac{r_+}{s_0})'}$. Note that the latter
  condition can be written as
  $w^{r_0}\in A_{\frac{r_0}{p_-}}\cap RH_{(\frac{p_+}{r_0})'}$ with
  $p_-=r_-$ and
  $\frac1{p_+}=\frac{1}{r_0}-\frac1{s_0}+\frac1{r_+}$, and in this
  case 
  $\frac1{s_0}-\frac1{r_0}+\frac1{p_+}=\frac1{r_+}\ge 0$, so the
  hypotheses of Theorem~\ref{thm:off-diag-limited} hold.

\medskip

A restricted range, off-diagonal extrapolation theorem
has previously appeared in the literature.
Duoandikoetxea~\cite[Theorem~5.1]{MR2754896} proved that if for some
$1\leq p_0<\infty$ and $0<q_0,\,r_0<\infty$, and all weights $w\in
A_{p_0,r_0}$ (note that unlike in the classical definition of this
class he does not require $p_0\leq q_0$), if
\eqref{eqn:off-diag-limited1} holds, then for all $1<p<\infty$ and
$0<q,\,r<\infty$ such that
$\frac{1}{p_0}-\frac{1}{p}=\frac{1}{q_0}-\frac{1}{q}=\frac{1}{r_0}-\frac{1}{r}$,
and all weights $w\in A_{p,r}$, \eqref{eqn:off-diag-limited2} holds.

This result is contained in Theorem~\ref{thm:off-diag-limited} in the
particular case when  $r_0\ge \min\{p_0,q_0\}$ if we take $p_-=1$ and $p_+=\big(\frac1{p_0}-\frac1{r_0}\big)^{-1}$. In this case,  (because $r_0\geq p_0$)
$w\in A_{p_0, r_0}$ if and only if
$w^{p_0}\in A_{p_0}\cap RH_{\frac{r_0}{p_0}}= A_{\frac{p_0}{p_-}} \cap
RH_{\big(\frac{p_+}{p_0}\big)'}$.   Moreover, in this scenario 
$\frac1{q_0}-\frac1{p_0}+\frac1{p_+}\ge 0$  since $r_0\ge q_0$. 

Despite this overlap, our results are different.  We eliminate
the restriction $p_0,\,p>1$ as we can take $0\leq p_-<1$.  Moreover,
for a value of $p_-\neq 1$, it is not clear whether our result can be
gotten from his by rescaling.  On the other hand, we cannot recapture
his result for values of $r_0<\min\{p_0,q_0\}$.  

Finally, in light of Remark~\ref{remark:rest-expon:0}, we note that
\cite[Theorem~5.1]{MR2754896} allows for weights $w^{q_0}$ or
$w^{p_0}$ that may not be locally integrable unless one assumes
$r_0\ge \min\{p_0,q_0\}$. For example, if we fix
$0<r_0<\min\{p_0,q_0\}$ and let
$w(x)=|x|^{-\frac{n}{\min\{p_0,q_0\}}}$, then it is easy to see that
$w^{r_0}\in A_1$ and so $w\in A_{p_0,r_0}$, but neither $w^{p_0}$ nor
$w^{q_0}$ is locally integrable (and so the characteristic function of
the unit ball centered at $0$ does not belong to $L^{p_0}(w^{p_0})$ or
to $L^{q_0}(w^{q_0})$). In light of this, we believe the condition
$r_0\ge \min\{p_0,q_0\}$ is not unduly restrictive.

\medskip

\subsection*{Applications}
To demonstrate the power of our multilinear extrapolation theorem, we
use Theorem~\ref{thm:extrapol-multi} to prove results for the bilinear
Hilbert transform and for multilinear Calder\'on-Zygmund operators.
We first consider the bilinear Hilbert transform, which is defined by
\[ BH(f_1,f_2)(x) = 
\text{p.v.} \int_{\re} f_1(x-t)g(x+t)\frac{dt}{t}.  \]
The problem of finding bilinear $L^p$ estimates for this operator was
first raised by Calder\'on in connection with the Cauchy integral
problem (though it was apparently not published
until~\cite{MR734178}).  Lacey and Thiele~\cite{MR1491450,MR1689336} showed
that for $1<p_1,\,p_2\le\infty$,
$\frac{1}{p}=\frac{1}{p_1}+\frac{1}{p_2}<\frac{3}{2}$,
\[ \|BH(f_1,f_2)\|_{L^p} 
\leq C\|f_1\|_{L^{p_1}}\|f_2\|_{L^{p_2}}. \]

The problem of weighted norm inequalities for the bilinear Hilbert
transform has been raised by a number of authors:
see~\cite{MR2997585,MR3000985,grafakos-martell04,MR3130311}.  The
first such results were recently obtained by
Culiuc, di Plinio and Ou~\cite{Culiuc:2016wr}.  

\begin{theorem} \label{thm:bht-cdo}
Given $1<p_1,\,p_2<\infty$, define $p$ by
$\frac{1}{p}=\frac{1}{p_1}+\frac{1}{p_2}$ and assume  that
$p>1$.  For $i=1,\,2$, let $w_i$ be such that $w_i^{2p_i} \in
A_{p_i}$, and define $w=w_1w_2$.  
Then 
\begin{equation} \label{eqn:bht0}
\|BH(f_1,f_2)\|_{L^p(w^p)} \leq
    C\|f_1\|_{L^{p_1}(w_1^{p_1})}\|f_2\|_{L^{p_2}(w_2^{p_2})}, 
  \end{equation}
where $C=C(p_i,[w_i^{2p_i}]_{A_{p_i}})$.
\end{theorem}

If we apply   Theorem~\ref{thm:extrapol-multi}, we can extend
Theorem~\ref{thm:bht-cdo} to a larger collection of weights and exponents.  In
particular, we can remove the restriction that $p>1$, replacing it
with $p>\frac{2}{3}$, the same threshold that appears in the unweighted
theory.

\begin{theorem} \label{thm:bht-DCU-JMM}
Given arbitrary $1<p_1,\,p_2<\infty$, define $\frac{1}{p}=\frac{1}{p_1}+\frac{1}{p_2}$ and assume that $p>1$. For every $i=1,2$, let $r_i^-=\frac{2 p_i}{1+p_i}<q_i<2p_i=r_i^+$. Then, for all 
$w_i^{q_i}\in A_{\frac{q_i}{r_i^-}}\cap RH_{\big(\frac{r_i^+}{q_i}\big)'}$ ---or, equivalently, $w_i^{2r_i}\in A_{r_i}$ for $r_i=\big(\frac2{q_i}-\frac1{p_i}\big)^{-1}$---
if we write $w=w_1w_2$ and $\frac1q=\frac1{q_1}+\frac1{q_2}$, we have that
\begin{equation} \label{eqn:bht1-chema}
\|BH(f_1,f_2)\|_{L^{q}(w^q)}
\leq C\|f_1\|_{L^{q_1}(w_1^{q_1})}
\|f_2\|_{L^{q_2}(w_2^{q_2})}.
\end{equation}

In particular, given arbitrary $1<q_1,\,q_2<\infty$ so that
$q>\frac23$ where $\frac{1}{q}=\frac{1}{q_1}+\frac{1}{q_2}$, there exist
values $1<p_1,\,p_2<\infty$ such that
$\frac{1}{p}=\frac{1}{p_1}+\frac{1}{p_2}<1$, in such a way that if we
set $r_i^-=\frac{2 p_i}{1+p_i}$, $r_i^+=2p_i$ then $r_i^-<q_i<r_i^+$,
and for all weights $w_i$ with
$w_i^{q_i}\in A_{\frac{q_i}{r_i^-}}\cap
RH_{\big(\frac{r_i^+}{q_i}\big)'}$ (or, equivalently,
$w_i^{2r_i}\in A_{r_i}$ for
$r_i=\big(\frac2{q_i}-\frac1{p_i}\big)^{-1}$) and $w=w_1w_2$,
\begin{equation} \label{eqn:bht1}
\|BH(f_1,f_2)\|_{L^{q}(w^q)}
\leq C\|f_1\|_{L^{q_1}(w_1^{q_1})}
\|f_2\|_{L^{q_2}(w_2^{q_2})}.
\end{equation}
\end{theorem}

\begin{remark}
  We can state Theorem~\ref{thm:bht-DCU-JMM} in a different but
  equivalent form.  For instance, in the second part of that result, if  we let $v_i=w_i^{q_i}$, then our hypothesis
  becomes
  $v_i \in A_{\frac{q_i}{r_i^-}}\cap RH_{\big(\frac{r_i^+}{q_i}\big)'}$, and
  the conclusion is that
\[ BH : L^{q_1}(v_1) \times L^{q_2}(v_2) 
\longrightarrow L^q(v_1^{\frac{q}{q_1}}v_w^{\frac{q}{q_2}}). \]
In~\cite{Culiuc:2016wr}, for instance, Theorem~\ref{thm:bht-cdo} is
stated in this form.  We chose the form that we did because it seems
more natural when working with off-diagonal inequalities.
\end{remark}

\begin{remark}
In~\cite{Culiuc:2016wr} the authors actually proved
Theorem~\ref{thm:bht-cdo} for a more general
family of bilinear multiplier operators introduced by Muscalu, Tao and
Thiele~\cite{Muscalu:2002wx}.  Theorem~\ref{thm:bht-DCU-JMM}
immediately extends to these operators.  We refer the interested
reader to these papers for precise definitions.
This extension actually shows that that the bound $p>1$ in
Theorem~\ref{thm:bht-cdo} and the bound $p>\frac{2}{3}$ in
Theorem~\ref{thm:bht-DCU-JMM} are natural and in some sense the best
possible.   In~\cite[Theorem~2.14]{MR1745019}, Lacey gave an example
of an operator which does not satisfy a bilinear estimate when $p<2/3$;
in~\cite[Remark~1.2]{Culiuc:2016wr} the authors show that
Theorem~\ref{thm:bht-cdo} applies to this operator.   Hence, if Theorem~\ref{thm:bht-cdo} could
be extended to include the case $p<1$, we would get weighted estimates
for this operator.
But by extrapolation, these would yield inequalities below the
threshold $q=\frac{2}{3}$.  Indeed, we could apply the
  first part of Theorem \ref{thm:bht-DCU-JMM} with those fixed
  exponents $\frac1p=\frac1{p_1}+\frac1{p_2}>1$ and $w_1=w_2\equiv 1$
  to obtain that this operator maps
 $L^{q_1}\times L^{q_2}$ into $L^q$ for every
  $r_i^-=\frac{2 p_i}{1+p_i}<q_i<2p_i=r_i^+$ and
  $\frac1q=\frac1{q_1}+\frac1{q_2}$. If we fix
  $0<\epsilon<\min\{\frac12(\frac1{p}-1),\frac1{p_i'}\}$ and let
  $\frac{1}{q_i}:=\frac12(\frac1{p_i}+1-\epsilon)$, we would have that
  $r_i^-<q_i<r_i^+$ and
$$
\frac1q
=
\frac1{q_1}+\frac1{q_2}
=
\frac1{2p}+1-\epsilon
>\frac32.
$$
\end{remark}

\bigskip

Given $q_1,\,q_2$, as part of the proof of Theorem~\ref{thm:bht-DCU-JMM} we construct the
parameters $r_i^-,r_i^+$ needed to define the weight classes.   Thus,
while we show that such weights exist, it is not clear from the
statement of the theorem what weights are possible.  To illustrate the different kinds of weight
conditions we get, we give some special classes of weights, and in particular we give a family of power weights. 

\begin{cor} \label{cor:bht-A1}
Given $1<q_1,\,q_2<\infty$, define $q$ by
$\frac{1}{q}=\frac{1}{q_1}+\frac{1}{q_2}$, and assume further that
$q>\frac{2}{3}$.  Then,  
\begin{equation} \label{eqn:bht1-bis}
\|BH(f_1,f_2)\|_{L^{q}(w^q)}
\leq C\|f_1\|_{L^{q_1}(w_1^{q_1})}
\|f_2\|_{L^{q_2}(w_2^{q_2})}
\end{equation}
holds for all $w_i^{q_i}\in  A_{\max\{1,\frac{q_i}2\}}\cap RH_{\max\{1,\frac2{q_i}\}}$ 
and $w=w_1w_2$. 
In particular, 
\begin{equation}
BH : L^{q_1}(|x|^{-a}) \times L^{q_2}(|x|^{-a}) 
\longrightarrow L^q(|x|^{-a}),
\label{eq:BH-power}
\end{equation}
if $a=0$ or if 
\begin{equation}\label{eq:values-a}
1-\min\left\{ 
\max\left\{1,\frac{q_1}{2}\right\},\max\left\{1,\frac{q_2}{2}\right\}\right\}
< 
a
<\min\left\{1, \frac{q_1}{2},\frac{q_2}{2}\right\}.
\end{equation}
As a result, \eqref{eq:BH-power} holds for all $0\le a<\frac12$.
\end{cor}

\begin{remark}
By Corollary~\ref{cor:bht-A1} we get  weighted estimates for the bilinear Hilbert transform
  in exactly the same range where the unweighted estimates are known to hold.
  (Note that when $a=0$ we recover the unweighted case.) 
  Rather than taking equal weights in \eqref{eq:BH-power}, we can also
  give this inequality for more general power weights of the form
  $w_i=|x|^{-a_i/q_i}$; details are left to the interested reader.
\end{remark}

\begin{remark}\label{remark:A1} 
  As a consequence of Corollary~\ref{cor:bht-A1} we see that even in
  the range of exponents covered by Theorem \ref{thm:bht-cdo}
  from~\cite{Culiuc:2016wr}, we get a larger class of weights.  Fix
  $1<q_1,q_2<\infty$ and assume that
  $\frac{1}{q}=\frac{1}{q_1}+\frac{1}{q_2}<1$ . First, it is easy to
  show (see Lemma \ref{lemma:cjn} below) that $w_i^{2q_i}\in A_{q_i}$
  if an only if $w_i^{q_i}\in A_{\frac{1+q_i}{2}}\cap RH_2$. Hence, if
  we further assume that $w_i^{q_i}\in A_1$ this condition becomes
  $w_i^{q_i}\in A_1\cap RH_2$ or, equivalently, (see Lemma
  \ref{lemma:cjn} below) $w_i^{2q_i}\in A_1$. Hence, as a corollary of
  Theorem \ref{thm:bht-cdo} we get that
  $BH:L^{q_1}(w_1^{q_2})\times L^{q_2}(w_2^{q_2})\longrightarrow
  L^q(w^q)$ for all $w_i^{2q_i}\in A_1$. But by
  Corollary~\ref{cor:bht-A1}, again assuming that
  $w_i^{q_i}\in A_1$, we can allow
  $w^{q_i}\in A_1\cap RH_{\max\{1,\frac2{q_i}\}}$, or equivalently,
  $w_i^{\max\{2,q_i\}}\in A_1$ which is weaker than
  $w_i^{2q_i}\in A_1$ since $\max\{2,q_i\}< 2q_i$.
	
  Further, when $1<q_i\le 2$, Corollary~\ref{cor:bht-A1} gives the
  class of weights $w_i^{q_i}\in A_1\cap RH_{\frac2{q_i}}$. To compare
  this with Theorem \ref{thm:bht-cdo} from~\cite{Culiuc:2016wr} note
  that their condition is, as explained above,
  $w_i^{q_i}\in A_{\frac{1+q_i}{2}}\cap RH_2$ and hence we can weaken
  $w^{q_i}\in RH_2$ to $w^{q_i}\in RH_{\frac2{q_i}}$ at the cost of
  assuming that $w^{q_i}\in A_1$ . Alternatively, if $q_i\ge 2$, our
  condition becomes $w^{q_i}\in A_{\frac{q_i}2}$, which removes any
  reverse H\"older condition for $w^{q_i}$ at the cost of assuming
  that $w^{q_i}\in A_{\frac{q_i}2}\subset A_{\frac{1+q_i}{2}}$.
\end{remark}

\medskip

We can
also prove vector-valued inequalities for the bilinear Hilbert
transform for the same weighted Lebesgue spaces as
in the scalar inequality.   Even in the unweighted case, vector-valued
inequalities were an
open question until recently.  Benea and
Muscalu~\cite{MR3599522, Benea:2016th} (see also \cite{MR3124934,
  MR3291796} for earlier results) proved that given
$1<s_1,\,s_2\le \infty$ and $s$ such that
$\frac{1}{s}=\frac{1}{s_1}+\frac{1}{s_2}$ and $s>\frac{2}{3}$, then
there exist $q_1,\,q_2,\,q$ such that
\[ \bigg\| \bigg(\sum_k |BH(f_k,g_k)|^s\bigg)^{\frac{1}{s}}
\bigg\|_{q} 
\leq 
C\bigg\| \bigg(\sum_k |f_k|^{s_1}\bigg)^{\frac{1}{s_1}}
\bigg\|_{q_1} 
\bigg\| \bigg(\sum_k |g_k|^{s_2}\bigg)^{\frac{1}{s_2}}
\bigg\|_{q_2} \]
where $1<q_1,\,q_2\le\infty$, $\frac{1}{q}=\frac{1}{q_1}+\frac{1}{q_2}$,
and, depending on the values of the $s_i$, there are additional
restrictions on the possible values of the $q_i$.
(See~\cite[Theorem~5]{Benea:2016th} for a precise statement or \eqref{eq:new-cond-vv:1} below.)  An alternative proof of these estimates when $s>1$ is given in~\cite{Culiuc:2016wr}.

By using the formalism of extrapolation pairs, vector-valued
inequalities are an immediate consequence of extrapolation.  Hence, as
a consequence of Theorem~\ref{thm:bht-DCU-JMM} we get the following
generalization of the results
in~\cite{MR3599522,Benea:2016th,Culiuc:2016wr}. We note that for some triples $s_1,\,s_2,\,s$ our
method does not let us recover the full range of spaces gotten
in~\cite{MR3599522, Benea:2016th} but we do get
weighted estimates in our range.

\begin{theorem} \label{thm:bht-vector} 
Given arbitrary $1<p_1,\,p_2<\infty$, define $\frac{1}{p}=\frac{1}{p_1}+\frac{1}{p_2}$ and assume that $p>1$. For every $i=1,2$, let $r_i^-=\frac{2 p_i}{1+p_i}<q_i,s_i<2p_i=r_i^+$. Then,  for all 
$w_i^{q_i}\in A_{\frac{q_i}{r_i^-}}\cap RH_{\big(\frac{r_i^+}{q_i}\big)'}$ ---or, equivalently, $w_i^{2r_i}\in A_{r_i}$ for $r_i=\big(\frac2{q_i}-\frac1{p_i}\big)^{-1}$---
if we write $w=w_1w_2$, $\frac1q=\frac1{q_1}+\frac1{q_2}$ and  $\frac1s=\frac1{s_1}+\frac1{s_2}$, there holds
\begin{multline} \label{eqn:bht-vector1}
\bigg\| \bigg(\sum_k |BH(f_k,g_k)|^s\bigg)^{\frac{1}{s}}
\bigg\|_{L^q(w^q)} \\
\leq 
C\bigg\| \bigg(\sum_k |f_k|^{s_1}\bigg)^{\frac{1}{s_1}}
\bigg\|_{L^{q_1}(w_1^{q_1})} 
\bigg\| \bigg(\sum_k |g_k|^{s_2}\bigg)^{\frac{1}{s_2}}
\bigg\|_{L^{q_2}(w_2^{q_2})}. 
\end{multline}
In particular, for every $1<s_1,\,s_2<\infty$ such that $\frac{1}{s}=\frac{1}{s_1}+\frac{1}{s_2}<\frac32$, and for every $1<q_1, q_2<\infty$ such that
  $\frac{1}{q}=\frac{1}{q_1}+\frac{1}{q_2}<\frac{3}{2}$, if
	\begin{equation}	\label{eq:restriction-exponents-vv}
	\left|\frac1{s_1}-\frac1{q_1}\right|<\frac12,
	\qquad\quad
	\left|\frac1{s_2}-\frac1{q_2}\right|<\frac12,
	\qquad\text{and}\qquad
	\sum_{i=1}^2 \max\left\{\frac1{q_i},\frac1{s_i}\right\}<\frac32,
\end{equation}
 there are
values $1<p_1,\,p_2<\infty$ such that
$\frac{1}{p}=\frac{1}{p_1}+\frac{1}{p_2}<1$, in such a way that if we
set $r_i^-=\frac{2 p_i}{1+p_i}$, $r_i^+=2p_i$ then $r_i^-<q_i,s_i<r_i^+$, and hence
\eqref{eqn:bht-vector1} holds for all weights $w_i$ with
$w_i^{q_i}\in A_{\frac{q_i}{r_i^-}}\cap
RH_{\big(\frac{r_i^+}{q_i}\big)'}$ (or, equivalently,
$w_i^{2r_i}\in A_{r_i}$ for
$r_i=\big(\frac2{q_i}-\frac1{p_i}\big)^{-1}$) and $w=w_1w_2$. 
\end{theorem}

\medskip

\begin{remark}
  Theorem~\ref{thm:bht-vector} contains the vector-valued inequalities that
  follow immediately from our extrapolation result applied to the
  weighted norm inequalities obtained in~\cite{Culiuc:2016wr}
  (cf. Therorem~\ref{thm:bht-cdo}). However, more general weighted
  estimates for the bilinear Hilbert transform are implicit in the
  arguments of~\cite{Culiuc:2016wr}. These in turn produce
  vector-valued inequalities in a wider range of exponents. We shall
  elaborate on this in Section~\ref{section:generalize} below.
\end{remark}

\begin{remark} \label{remark:iteration}
In~\cite[Proposition~10]{MR3599522} the authors also prove iterated
vector-valued inequalities of the form
\begin{multline*}
 \bigg\| \bigg(\sum_j 
\bigg(\sum_k |BH(f_{jk},g_{jk})|^s\bigg)^{\frac{t}{s}}
\bigg)^{\frac{1}{t}}
\bigg\|_{p} \\
\leq 
\bigg\| \bigg(\sum_j 
\bigg(\sum_k |f_{jk}|^{s_1}\bigg)^{\frac{t_1}{s_1}}
\bigg)^{\frac{1}{t_1}}
\bigg\|_{p_1} 
\bigg\| \bigg(\sum_j 
\bigg(\sum_k |g_{jk}|^{s_2}\bigg)^{\frac{t_2}{s_2}}
\bigg)^{\frac{1}{t_2}}
\bigg\|_{p_2}, 
\end{multline*}
again with restrictions on the possible values of the $p_i$ depending
on the $s_i$ and $t_i$.  We can easily prove some of these
inequalities \ by extrapolation; moreover, we can also prove prove
weighted versions.  After the proof of Theorem~\ref{thm:bht-vector} we
sketch how this is done.   Here we note in passing that iterated vector-valued
inequalities have recently appeared in another setting:  see~\cite{2017arXiv170507792A}.
\end{remark}

\medskip 

As we did with the scalar inequalities we give some specific examples
of classes
of weights for which the bilinear Hilbert transform satisfies weighted
vector-valued inequalities.

\begin{cor}\label{corol:BHT-vv-power}
 Given $1<s_1,\,s_2<\infty$ such that $\frac{1}{s}=\frac{1}{s_1}+\frac{1}{s_2}<\frac32$, and $1<q_1, q_2<\infty$ such that   $\frac{1}{q}=\frac{1}{q_1}+\frac{1}{q_2}<\frac{3}{2}$, if
	\begin{equation}	\label{eq:restriction-exponents-vv-power-corol}
	\left|\frac1{s_1}-\frac1{q_1}\right|<\frac12,
	\qquad\quad
	\left|\frac1{s_2}-\frac1{q_2}\right|<\frac12,
	\qquad\text{and}\qquad
	\sum_{i=1}^2 \max\left\{\frac1{q_i},\frac1{s_i}\right\}<\frac32,
\end{equation}
 then \eqref{eqn:bht-vector1} holds for all $w_i^{q_i}\in A_{\max\{1,\frac{q_i}2,\frac{q_i}{s_i}\}}\cap RH_{\max\{1,\frac2{q_i}, [1-q_i(\frac1{s_i}-\frac12)]^{-1}\}}$. In particular, 
\begin{multline} \label{eqn:bht-vector1:pwerw-wts}
\bigg\| \bigg(\sum_k |BH(f_k,g_k)|^s\bigg)^{\frac{1}{s}}
\bigg\|_{L^q(|x|^{-a})} \\
\leq 
C\bigg\| \bigg(\sum_k |f_k|^{s_1}\bigg)^{\frac{1}{s_1}}
\bigg\|_{L^{q_1}(|x|^{-a} )} 
\bigg\| \bigg(\sum_k |g_k|^{s_2}\bigg)^{\frac{1}{s_2}}
\bigg\|_{L^{q_2}(|x|^{-a})}. 
\end{multline}
holds if $a\in\{0\}\cup(a_-,a_+)$ where
\begin{align}
\label{eq:rest-a-v-v}
a_-&=1-\min \left\{ 
\max\left\{1,\frac{q_1}2, \frac{q_1}{s_1}\right\},
\max\left\{1,\frac{q_2}2, \frac{q_2}{s_2}\right\}
\right\}
\\[6pt] \nonumber
a_+&=
\min \left\{
1,
\frac{q_1}2,
\frac{q_2}2, 
1-q_1\left(\frac1{s_1}-\frac12\right), 
1-q_2\left(\frac1{s_2}-\frac12\right)\right\}
\end{align}
\end{cor}

\medskip

\begin{remark} 
  The conditions in~\eqref{eq:restriction-exponents-vv-power-corol}
  guarantee that $a_-\le 0<a_+$, hence the set $\{0\}\cup(a_-,a_+)$
  defines a non-empty interval. On the other hand, this interval can
  be arbitrarily small. For instance, take $q_1=s_1=2$, $q_2=2$,
  $s_2=t$ with $1<t<2$. Then
  \eqref{eq:restriction-exponents-vv-power-corol} is satisfied and we
  have that $a_-=0$ and $a_+=2(1-\frac1t)$. Thus,
  $\{0\}\cup(a_-,a_+)=[0,a_+)$ and $a_+\to 0$ as $t\to 1^+$: that is,
  in the limit we just get the Lebesgue measure. Notice, however, that
  in the context of the first part of
  Corollary~\ref{corol:BHT-vv-power}, as $t\to 1^+$, the conditions on
  the weights become $w_1^2\in A_1$ and $w_2^2\in A_2\cap
  RH_\infty$. Hence, we can take $w_1(x)=|x|^{-\frac{a_1}{q_1}}$ and
  $w_2(x)=|x|^{-\frac{a_2}{q_2}}$ with $a_1\ge 0$ and $-1<a_2\le 0$.
  (Of course if $a_1=a_2=a$, then $a=0$ as observed above.)
\end{remark}

\medskip

As a final application we use extrapolation to prove Marcinkiewicz-Zygmund
inequalities for multilinear Calder\'on-Zygmund operators.   Weighted
norm inequalities for these operators have been considered by several
authors:  we refer the reader to~\cite{grafakos-martell04,MR2483720}
for precise definitions of these operators and weighted norm
inequalities for them.  Very recently, Carando, Mazzitelli and
Ombrosi~\cite{Carando:2016vm} proved the following weighted
Marcinkiewicz-Zygmund inequalities.

\begin{theorem} \label{theorem:CZ-MZ}
For $m\geq 1$, let $T$ be an $m$-linear Calder\'on-Zygmund operator.
Given $1<q_1,\ldots,q_m<\infty$, $q$ such that $\frac{1}{q}=\sum
\frac{1}{q_i}$, and weights $w_i$ such that
$w_i^{q_i} \in A_{q_i}$, 
\begin{equation} \label{eqn:CZ-MZ0}
\bigg\| \bigg(\sum_{k_1,\ldots,k_m} 
|T(f_{k_1}^1,\ldots,f_{k_m}^m)|^2\bigg)^{\frac{1}{2}}
\bigg\|_{L^q(w^q)} 
\leq C\prod_{i=1}^m
\bigg\| \bigg( \sum_{k_i} |f_{k_i}^i|^2\bigg)^{\frac{1}{2}}
\big\|_{L^{q_i}(w_i^{q_i})}, 
\end{equation}
where $w=w_1w_2$.
If $1<r<2$ and if we further assume $1<q_i<r$, then again for all
weights $w_i$ such that $w_i^{q_i} \in A_{q_i}$,
\begin{equation} \label{eqn:CZ-MZ1}
\bigg\| \bigg(\sum_{k_1,\ldots,k_m} 
|T(f_{k_1}^1,\ldots,f_{k_m}^m)|^r\bigg)^{\frac{1}{r}}
\bigg\|_{L^q(w^q)} 
\leq C\prod_{i=1}^m
\bigg\| \bigg( \sum_{k_i} |f_{k_i}^i|^r\bigg)^{\frac{1}{r}}
\big\|_{L^{q_i}(w_i^{q_i})},
\end{equation}
where $w=w_1w_2$.
\end{theorem}

By using extrapolation we can prove that inequality~\eqref{eqn:CZ-MZ1}
holds for $1<r<2$ with the same family of exponents as
in~\eqref{eqn:CZ-MZ0} for $r=2$.  

\begin{theorem} \label{thm:CZ-MZ-new}
For $m\geq 1$, let $T$ be an $m$-linear Calder\'on-Zygmund operator.
Given $1<r<2$,  $1<q_1,\ldots,q_m<\infty$, $q$ such that $\frac{1}{q}=\sum
\frac{1}{q_i}$, and weights $w_i$ such that
$w_i^{q_i} \in A_{q_i}$, then inequality~\eqref{eqn:CZ-MZ1} holds.
\end{theorem}

\begin{remark}
In~\cite{Carando:2016vm} the authors actually prove that
Theorem~\ref{theorem:CZ-MZ} holds for weights in the larger class
$A_{\vec{p}}$ introduced in~\cite{MR2483720}.  However, it is not
known whether multilinear extrapolation holds for these weights.   We
also do not know if Theorem~\ref{thm:CZ-MZ-new} can be extended to
this larger family of weights.
\end{remark}

\medskip

The remainder of this paper is organized as follows.  In
Section~\ref{section:prelim} we gather some definitions and basic
results about weights.  In Section~\ref{section:extrapol} we prove all
of our extrapolation results. In
Section~\ref{section:BHT} we give the proofs of all of the
applications.  Finally, in Section~\ref{section:generalize} we discuss
some results that are implicit in~\cite{Culiuc:2016wr} and that can be
used to get more general vector-valued inequalities for the bilinear
Hilbert transform.

Throughout this paper $n$ will denote the dimension of the underlying
space, $\rn$.  A constant $C$ may depend on the dimension $n$, the
underlying parameters $p_-,p_+,p,\ldots$, and the $A_p$ and $RH_s$
constants of the associated weights.  It will not depend on the
specific weight.  The value of a constant $C$ may change from line to
line.  Throughout we will  use the conventions that
$\frac{1}{\infty}=0$, $\frac{1}{0}=\infty$, and $1'=\infty$ and
$\infty'=1$. 

\section{Preliminaries}
\label{section:prelim}

In this section we give the basic properties of weights that we will
need below.  For proofs and further information,
see~\cite{duoandikoetxea01,grafakos08b}.   By a weight we mean a
non-negative  function $v$ such that
$0<v(x)<\infty$ a.e.   For $1<p<\infty$, we say $v\in
A_p$ if 
\[ [v]_{A_p} = \sup_Q \avgint_Q v\,dx 
\left(\avgint_Q v^{1-p'}\,dx\right)^{p-1} < \infty, \]
where the supremum is taken over all cubes $Q\subset \rn$ with sides
parallel to the coordinate axes and $\avgint_Q v\,dx= |Q|^{-1}\int_Q
v\,dx$.   The quantity $[v]_{A_p}$ is called
the $A_p$ constant of $v$.  Note that it follows at once from this
definition that if $v\in A_p$, then $v^{1-p'}\in A_{p'}$.  When $p=1$
we say $v\in A_1$ if
\[  [v]_{A_1} = \sup_Q \avgint_Q w(y)\,dy \esssup_{x\in Q} w(x)^{-1}<
  \infty. \]
The $A_p$ classes are properly nested:  for $1<p<q$, $A_1\subsetneq
A_p \subsetneq A_q$.  
We denote the union of all the $A_p$ classes, $1\leq p<\infty$, by
$A_\infty$.  

Given $1<s<\infty$, we say that a weight $v$ satisfies the reverse
H\"older inequality with exponent $s$, denoted $w \in RH_s$ if
\[ [v]_{RH_s} = \sup_Q \left(\avgint_Q v^s\,dx\right)^{\frac{1}{s}} \left(\avgint_Q v\,dx\right)^{-1}
< \infty. \]
When $s=\infty$ we say $v\in RH_\infty$ if
\[ [ v]_{RH_\infty} = \sup_Q \esssup_{x\in Q} w(x) 
\left(\avgint_Q v\,dx\right)^{-1} < \infty.  \]
The reverse H\"older classes are also properly nested:  if
$s<t<\infty$, then $RH_\infty \subsetneq RH_t \subsetneq RH_s$. Define $RH_1$ to be the union of all the
$RH_s$ classes, $1<s\leq \infty$.    We have that
$RH_1=A_\infty$.  A given 
$v$ is in $RH_s$ for some $s>1$ if and only if there exists $p>1$ such
that $v\in A_p$.   Equivalently,  if $v\in A_\infty$, there exists
$1\leq p<\infty$ and $1<s\leq \infty$ such that 
$v\in A_p \cap RH_s$.  

The $A_p$ and $RH_s$ classes satisfy openness properties:  given $v
\in A_p$, $1<p<\infty$, then there exists $\epsilon>0$
depending only on $[v]_{A_p}$, $p$ and $n$, such that $v\in A_{p-\epsilon}$;   also given $v
\in RH_s$, $1<s<\infty$, then there exists $\epsilon>0$
depending only on $[v]_{RH_s}$, $s$, and $n$, such that $v\in RH_{s+\epsilon}$.

The condition $v\in A_p \cap RH_s$ can be restated using the following
result.  The first part is from~\cite[Theorem~2.2]{MR1308005}; the
second is just gotten by the duality of $A_p$ weights.

\begin{lemma} \label{lemma:cjn}
Given $1\le p<\infty$, $1\le s<\infty$, the weight $v\in A_p \cap RH_s$ if and only if $v^s\in A_q$, where
$q= s(p-1)+1$, that is,
\begin{equation}
\sup_Q \left(\avgint_Q v^s\,dx\right)^\frac1s
\left(\avgint_Q v^{1-p'}\,dx\right)^{p-1} < \infty.
\label{eq:ApRHs}
\end{equation}
In this case also have that $v^{1-p'}\in A_{q'}$. 
\end{lemma}

We can also easily construct weights $v\in A_p \cap RH_s$.  The next result can be proved directly from the definitions of the weight classes; essentially the same argument is used to prove the easier half of the Jones factorization theorem.  See~\cite[Theorem~5.1]{MR1308005} or~\cite[Theorem~4.4]{2017arXiv170602620C}.

\begin{lemma} \label{lemma:rev-fac}
Given weights $v_1,\,v_2\in A_1$, then for all $1\le p<\infty$, $1<s \le\infty$,
\[ v = v_1^{\frac{1}{s}}v_2^{1-p} \in A_p\cap RH_s. \]
\end{lemma}

\section{Proofs of extrapolation results}
\label{section:extrapol}

Our proof is similar in spirit to the proofs of off-diagonal and
limited range extrapolation in~\cite[Theorems~3.23
and~3.31]{cruz-martell-perezBook}.  To better understand the heuristic
argument that underlies our proof, we refer the reader to the
discussion in~\cite[Section~4]{CruzUribe:2015km}.  We have split the
proof split into four cases.

\subsection*{Proof of Theorem~\ref{thm:off-diag-limited}. Case I: $\boldsymbol{p_->0}$ and $\boldsymbol{p_-<p_0<p_+}$}

Fix $p_-<p<p_+$ and $w$ such that
$w^p\in A_{\frac{p}{p_-}} \cap RH_{\big(\frac{p_+}{p}\big)'}$. Fix an
extrapolation pair $(f,g) \in \F$; we may assume that
$0<\|f\|_{L^q(w^q)},\, \|g\|_{L^p(w^p)}<\infty$.  For if
$\|f\|_{L^q(w^q)}=0$ or if $\|g\|_{L^p(w^p)}=\infty$, then
\eqref{eqn:off-diag-limited2} is trivially true. And if
$\|g\|_{L^p(w^p)}=0$, then \eqref{eqn:off-diag-limited1} implies that
$\|f\|_{L^{q_0}(w^{q_0})}=0$, and so $f=0$ a.e.~ and thus
$\|f\|_{L^{q}(w^{q})}=0$, which again gives us
\eqref{eqn:off-diag-limited2}.

\smallskip

We now fix some exponents based on our weight $w$.   By
Lemma~\ref{lemma:cjn} we have that $w^{p  \big(\frac{p_+}{p}\big)'}\in A_\tau$,
where 
\begin{equation} \label{eqn:tau}
 \tau = 
  \left(\frac{p_+}{p}\right)'\left(\frac{p}{p_-}-1\right)+1
= \frac{\frac{1}{p_-}-\frac{1}{p}}{\frac{1}{p}-\frac{1}{p_+}}+1
= \frac{\frac{1}{p_-}-\frac{1}{p_+}}{\frac{1}{p}-\frac{1}{p_+}}.
\end{equation}
For future reference we note that 
\begin{equation} \label{eqn:tau-prime}
\tau' =
  \frac{\frac{1}{p_-}-\frac{1}{p_+}}{\frac{1}{p_-}-\frac{1}{p}}. 
\end{equation} 
From Remark \ref{remark:rest-expon} we have that
\begin{equation}
\frac{1}{q_+}- \frac{1}{p_+}=\frac{1}{q_0}-\frac{1}{p_0}=\frac{1}{q}-\frac{1}{p}.
\label{eq:q+-right}
\end{equation}
Define the number $s$  by
\begin{equation} \label{eqn:s1}
 s = q_0 -
  \frac{q_0}{p_0}\frac{q}{\tau}\bigg(\frac{p_0}{p_-}-1\bigg)
= q_0q\bigg(\frac{1}{q}
-\frac{1}{\tau}\bigg(\frac{1}{p_-}-\frac{1}{p_0}\bigg)\bigg);
\end{equation}
we will explain our choice of $s$ below.  For later use, we prove that
$0<s< \min(q,q_0)$.  First, we have that $s>0$: by \eqref{eqn:tau},
the fact that $p_0<p_+$ and \eqref{eq:q+-right} we obtain
\begin{align*}
\frac{1}{q}-\frac{1}{\tau}\bigg(\frac{1}{p_-}-\frac{1}{p_0}\bigg)
 = \frac{1}{q}
  -\frac{\frac{1}{p}-\frac{1}{p_+}}{\frac{1}{p_-}-\frac{1}{p_+}}
\bigg(\frac{1}{p_-}-\frac{1}{p_0}\bigg)
 > \frac{1}{q} -\frac{1}{p}+\frac{1}{p_+} 
 = \frac{1}{q_+} \ge 0.
\end{align*}
To show that $s< \min(q,q_0)$, we claim 
\begin{equation} \label{eqn:s2}
 s = q-\frac{qq_0}{p_0}\frac{1}{\tau'}\frac{1}{\big(\frac{p_+}{p_0}\big)'}
= q_0q\bigg(\frac{1}{q_0}-\bigg(1-\frac{1}{\tau}\bigg)
\bigg(\frac{1}{p_0}-\frac{1}{p_+}\bigg)\bigg). 
\end{equation}
To see that this holds, we use the
fact that $\frac{1}{q}-\frac{1}{p}=\frac{1}{q_0}-\frac{1}{p_0}$:
\begin{align*}
 \frac{1}{q}
-\frac{1}{\tau}\bigg(\frac{1}{p_-}-\frac{1}{p_0}\bigg)
& = \frac{1}{q}
-\frac{1}{\tau}\bigg(\frac{1}{p_-}-\frac{1}{p_+}\bigg)
-\frac{1}{\tau}\bigg(\frac{1}{p_+}-\frac{1}{p_0}\bigg) \\
& = \frac{1}{q}-\frac{1}{p}+\frac{1}{p_+}
-\frac{1}{\tau}\bigg(\frac{1}{p_+}-\frac{1}{p_0}\bigg) \\
& = \frac{1}{q_0}-\frac{1}{p_0}+\frac{1}{p_+}
-\frac{1}{\tau}\bigg(\frac{1}{p_+}-\frac{1}{p_0}\bigg) \\
& = \frac{1}{q_0}
-\bigg(1-\frac{1}{\tau}\bigg)\bigg(\frac{1}{p_0}-\frac{1}{p_+}\bigg).
\end{align*}
It follows at once from \eqref{eqn:s1} and \eqref{eqn:s2} that 
$ s < \min(q,q_0)$.

\medskip

 We now prove our main estimate.  By rescaling and duality, we have that
\[ \|f\|_{L^q(w^q)}^s = \|f^s\|_{L^{\frac{q}{s}}(w^q)} 
= \int_\rn f^s h_2 w^q\,dx, \]
where $h_2$ is a non-negative function in $L^{(\frac{q}{s})'}(w^q)$ with
$\|h_2\|_{L^{(\frac{q}{s})'}(w^q)}=1$.  Now let $H_1$ and $H_2$ be
non-negative functions such that $0<H_1<\infty$ a.e., and $h_2\leq H_2$;
we will determine their exact values
 below.  Fix
$\alpha=\frac{s}{\big(\frac{q_0}{s}\big)'}$.  Then by H\"older's inequality,
\begin{multline} \label{eqn:splitting}
\int_\rn f^s h_2 w^q\,dx 
\leq \int_\rn f^s H_1^{-\alpha}H_1^\alpha H_2 w^q\,dx \\
\leq \left(\int_\rn f^{q_0}H_1^{-\alpha\frac{q_0}{s}}H_2
  w^q\,dx\right)^{\frac{s}{q}_0} 
\left( \int_\rn H_1^{\alpha\left(\frac{q_0}{s}\right)'}
  H_2w^q\,dx\right)^{1/\big(\frac{q_0}{s}\big)'} 
= I_1 ^{\frac{s}{q}_0} \times I_2^{1/\big(\frac{q_0}{s}\big)'}. 
\end{multline}

We first estimate $I_2$.  Assume that $H_1\in L^q(w^q)$ with
$\|H_1\|_{L^q(w^q)}\leq C_1<\infty$, and that $H_2 \in L^{\big(\frac{q}{s}\big)'}(w^q)$
with $\|H_2\|_{L^{\big(\frac{q}{s}\big)'}(w^q)}\leq C_2<\infty$.  Then again by
H\"older's inequality,
\begin{multline*}
I_2
\leq \left(\int_\rn H_1^{\alpha\left(\frac{q_0}{s}\right)'\frac{q}{s}}
w^q\,dx\right)^{\frac{s}{q}}
\left(\int_\rn
  H_2^{\big(\frac{q}{s}\big)'}w^q\,dx\right)^{1/\big(\frac{q}{s}\big)'} 
\\
\leq C_2 \left(\int_\rn H_1^q w^q\,dx\right)^{\frac{s}{q}}  \leq C_1^{s}C_2.
\end{multline*}

To estimate $I_1$ we want to apply~\eqref{eqn:off-diag-limited1}; to
do so we need to show that $I_1<\infty$.   Assume that  $f\leq
H_1\|f\|_{L^q(w^q)}$; then  we
have that 
\begin{multline*}
 I_1 \leq \|f\|_{L^q(w^q)}^{q_0}\int_\rn H_1^{q_0} H_1^{-\alpha\frac{q_0}{s}}H_2
  w^q\,dx \\
= \|f\|_{L^q(w^q)}^{q_0}\int_\rn H_1^s H_2
  w^q\,dx = \|f\|_{L^q(w^q)}^{q_0} \times I_2 < \infty. 
\end{multline*}

Define $\varphi=\big(\frac{q}{s}\big)'\frac{q_0}{p_0}$.
Then $\varphi>1$:  by \eqref{eqn:s2} we have that
\[ \frac{s}{q}\frac{p_0}{q_0} 
= \frac{p_0}{q_0}-\frac{1}{\tau'}\frac{1}{\big(\frac{p_+}{p_0}\big)'}, \]
and so
\[ \frac{1}{\varphi} = \frac{p_0}{q_0}\frac{1}{\big(\frac{q}{s}\big)'}=
\frac{p_0}{q_0}\bigg(1-\frac{s}{q}\bigg)=
\frac{1}{\tau'}\frac{1}{\big(\frac{p_+}{p_0}\big)'}< 1. \]
Now let $W^{q_0}= H_1^{-\alpha\frac{q_0}{s}}H_2 w^q$ and assume that
$W^{p_0}\in A_{\frac{p_0}{p_-}} \cap RH_{\big(\frac{p_+}{p_0}\big)'}$.  Since
$I_1$ is finite, $f\in L^{q_0}(W^{q_0})$.  Thus,
by~\eqref{eqn:off-diag-limited1} and H\"older's inequality,
\begin{align*}
I_1
& = \int_\rn f^{q_0} W^{q_0}\,dx  \\
& \leq C\left(\int_\rn g^{p_0} W^{p_0}\,dx\right)^{\frac{q_0}{p_0}} \\
& =  C\left(\int_\rn g^{p_0} H_1^{-\alpha \frac{p_0}{s}}
H_2^{\frac{p_0}{q_0}}w^{q\frac{p_0}{q_0}}w^{-q}w^q\,dx\right)^{\frac{q_0}{p_0}}
  \\ 
& \leq \left(\int_\rn g^{p_0\varphi'} H_1^{-\alpha \frac{p_0}{s}\varphi'}
w^{q(\frac{p_0}{q_0}-1)\varphi'}w^q\,dx\right)^{\frac{q_0}{\varphi' p_0}}
\left(\int_\rn  H_2^{\big(\frac{q}{s}\big)'}w^q\,dx\right)^{\frac{q_0}{\varphi p_0} }.
\end{align*}
The second integral on the last line is bounded by
$C_2^{\frac{q_0}{\varphi p_0}\big(\frac{q}{s}\big)'}=C_2$, so it remains
to show that the first integral is bounded by $\|g\|_{L^p(w^p)}^{q_0}$.   If we have that 
\[ g^{p_0\varphi'} H_1^{-\alpha \frac{p_0}{s}\varphi'}
w^{q(\frac{p_0}{q_0}-1)\varphi'} \leq H_1^q
\|g\|_{L^p(w^p)}^{p_0\varphi'}, \]
then the first integral would be bounded by
$\|H_1\|_{L^q(w^q)}^{\frac{q_0q}{\varphi' p_0}} \|g\|_{L^p(w^p)}^{q_0}\leq
C_1^{\frac{q_0q}{\varphi' p_0}}\|g\|_{L^p(w^p)}^{q_0}$.  This,
combined with inequality~\eqref{eqn:splitting} would yield
inequality~\eqref{eqn:off-diag-limited2} and the proof would be complete.

\medskip

Therefore, to complete the proof we need to show that we can construct
non-negative functions $H_1$ and $H_2$ such that 
\begin{gather}
\label{eqn:H1-norm}
\|H_1\|_{L^q(w^q)} \leq C_1, \\
\label{eqn:H1-pt}
g^{p_0\varphi'} H_1^{-\alpha \frac{p_0}{s}\varphi'}
w^{q(\frac{p_0}{q_0}-1)\varphi'} \leq H_1^q \|g\|_{L^p(w^p)}^{p_0\varphi'}, \\
\label{eqn:H1-f}
0<H_1<\infty, \quad f\leq H_1\|f\|_{L^q(w^q)},
\\
\label{eqn:H2-norm}
\|H_2\|_{L^{(\frac{q}{s})'}(w^q)} \leq C_2, \\
\label{eqn:H2-pt}
h_2 \leq H_2;
\end{gather}
and such that  the weight 
$W= H_1^{-\frac{\alpha}{s}}H_2^{\frac{1}{q_0}}  w^{\frac{q}{q_0}}$
satisfies
\begin{equation} \label{eqn:Ap}
W^{p_0}=
H_1^{-\frac{\alpha p_0}{s}}H_2^{\frac{p_0}{q_0}}
w^{\frac{qp_0}{q_0}}\in A_{\frac{p_0}{p_-}}\cap RH_{\big(\frac{p_+}{p_0}\big)'}. 
\end{equation}

\medskip

We will first prove that \eqref{eqn:H1-norm}, \eqref{eqn:H1-pt} and \eqref{eqn:H1-f} hold.
Since $\alpha\frac{p_0}{s}=\frac{p_0}{q_0}(q_0-s)$, 
\eqref{eqn:H1-pt} is equivalent to
\begin{equation} \label{eqn:H1-pt2}
 g^{p_0}w^{q(\frac{p_0}{q_0}-1)} 
\leq
H_1^{\frac{q}{\varphi'}+\frac{p_0}{q_0}(q_0-s)}\|g\|_{L^p(w^p)}^{p_0}. 
\end{equation}
Using the fact that
$\frac{1}{p}-\frac{1}{q}=\frac{1}{p_0}-\frac{1}{q_0}$, we have that
\begin{multline*}
\frac{q}{\varphi'}+\frac{p_0}{q_0}(q_0-s)
= q-\frac{p_0}{q_0}\frac{q}{\big(\frac{q}{s}\big)'}+\frac{p_0}{q_0}(q_0-s)
=q-\frac{p_0}{q_0}(q-s)+\frac{p_0}{q_0}(q_0-s) \\
=q-\frac{p_0}{q_0}(q-q_0) 
= q\bigg(1- p_0\bigg(\frac{1}{q_0}-\frac{1}{q}\bigg)\bigg)
= q\bigg(1- p_0\bigg(\frac{1}{p_0}-\frac{1}{p}\bigg)\bigg)
=q\frac{p_0}{p}.
\end{multline*}
Similarly, we have that 
\[ q\bigg(\frac{p_0}{q_0}-1\bigg)
= qp_0\bigg(\frac{1}{q_0}-\frac{1}{p_0}\bigg)
=qp_0\bigg(\frac{1}{q}-\frac{1}{p}\bigg)
= p_0 - q\frac{p_0}{p}
=
p_0\bigg(1-\frac{q}{p}\bigg)
. \]
Therefore, \eqref{eqn:H1-pt2} (and hence \eqref{eqn:H1-pt}) is equivalent to 
\begin{equation} \label{eqn:H1-pt3}
 g^{\frac{p}{q}} w^{\frac{p}{q}-1} \leq H_1 \|g\|_{L^p(w^p)}^{\frac{p}{q}}. 
\end{equation}

To construct a function $H_1$ that satisfies \eqref{eqn:H1-norm},
\eqref{eqn:H1-f}, and\eqref{eqn:H1-pt3}, 
we use the Rubio de Francia iteration algorithm.
As we noted above,  
$w^{p\big(\frac{p_+}{p}\big)'}\in A_\tau$, so  the maximal operator is bounded on 
$L^\tau(w^{p\big(\frac{p_+}{p}\big)'})$.  Hence,  for non-negative
$G\in L^\tau(w^{p(\frac{p_+}{p})'}) $ we can define the iteration
algorithm
\[ \R_1G = \sum_{k=0}^\infty 
\frac{M^k G}{2^k \|M\|_{L^\tau(w^{p(\frac{p_+}{p})'})}^k}.
 \]
Then we have that that $G\leq \R_1G$, $\R_1G\in A_1$, and 
$\|R_1 G\|_{L^\tau(w^{p(\frac{p_+}{p})'})} \leq 
2 \| G\|_{L^\tau(w^{p(\frac{p_+}{p})'})}$  
(cf.~\cite[Proof of Theorem~3.9]{cruz-martell-perezBook}).
Now define $\delta$ and $\epsilon$ by 
\[ \delta \tau =q, \qquad \epsilon
\tau = q-p\bigg(\frac{p_+}{p}\bigg)', \]
and let
\[ H_1 = \R_1(h_1^\delta
  w^\epsilon)^{\frac{1}{\delta}}w^{-\frac{\epsilon}{\delta}}, 
\qquad h_1 =
  \frac{f}{\|f\|_{L^q(w^q)}}+\frac{g^{\frac{p}{q}}w^{\frac{p}{q}-1}}{\|g\|_{L^p(w^p)}^{\frac{p}{q}}}. \]
Then
\[ \max\left(\frac{f}{\|f\|_{L^q(w^q)}},
\frac{g^{\frac{p}{q}}w^{\frac{p}{q}-1}}{\|g\|_{L^p(w^p)}^{\frac{p}{q}}}\right)
\leq h_1 \leq H_1, \]
and so both~\eqref{eqn:H1-f} and~\eqref{eqn:H1-pt3} hold.   Moreover,
\[ \|h_1\|_{L^q(w^q)} 
\leq 2^{1-\frac{1}{q}} \left(\int_\rn \frac{f^qw^q}{\|f\|_{L^q(w^q)}^q}
+ \frac{g^pw^p}{\|g\|_{L^p(w^p)}^p}\,dx\right)^{\frac{1}{q}} = 2, \]
and so
\begin{multline*}
\|H_1\|_{L^q(w^q)} 
=
\|\R_1(h_1^\delta
w^\epsilon)\|_{L^\tau(w^{p(\frac{p_+}{p})'})}^{\frac1\delta} \\
\leq 
2^{\frac1\delta}\|h_1^\delta w^\epsilon
\|_{L^\tau(w^{p(\frac{p_+}{p})'})}^{\frac1\delta}
= 2^{\frac1\delta}\|h_1\|_{L^q(w^q)}
\leq 2^{1+\frac1\delta} = C_1.
\end{multline*}
This gives us \eqref{eqn:H1-norm}. 

\medskip

The construction of $H_2$ and the proof of  \eqref{eqn:H2-norm} and
\eqref{eqn:H2-pt} are similar to the argument for $H_1$.   By Lemma \ref{lemma:cjn},  
if we set
\[ \sigma=p\bigg(\bigg(\frac{p}{p_-}\bigg)'-1\bigg), \]
 then $w^{-\sigma}\in
A_{\tau'}$ and so the maximal operator is bounded on
$L^{\tau'}(w^{-\sigma})$.  Hence, if we  define the Rubio de Francia iteration
algorithm for non-negative $F\in L^{\tau'}(w^{-\sigma})$ by 
\[ \R_2 F = \sum_{k=0}^\infty \frac{M^k F}{2^k
    \|M\|_{L^{\tau'}(w^{-\sigma})}^k},
			\]
then we have
that $F\leq \R_2F$, $\R_2F\in A_1$,  and $\|\R_2 F\|_{L^{\tau'}(w^{-\sigma})}
\leq 2\|F\|_{L^{\tau'}(w^{-\sigma})}$.  Define $\beta$ and $\gamma$ by
\[ \beta \tau' = \bigg(\frac{q}{s}\bigg)', \qquad 
\gamma \tau' = \sigma+q. \]
  If we now let
\[ H_2 = \R_2(h_2^\beta
  w^\gamma)^{\frac{1}{\beta}}w^{-\frac{\gamma}{\beta}}, \]
then we immediately get \eqref{eqn:H2-pt}.   Moreover, we have that
\begin{multline*}
\|H_2\|_{L^{(\frac{q}{s})'}(w^q)}
= \|\R_2(h_2^\beta
w^\gamma)\|_{L^{\tau'}(w^{-\sigma})}^{\frac1\beta} \\
\leq 2 ^{\frac1\beta}
\|h_2^\beta w^\gamma\|_{L^{\tau'}(w^{-\sigma})}^{\frac1\beta}
= 2 ^{\frac1\beta} \|h_2\|_{L^{(\frac{q}{s})'}(w^q)}
= 2 ^{\frac1\beta}  = C_2.
\end{multline*}
This gives us
\eqref{eqn:H2-norm}.

\medskip

Finally, we will show that \eqref{eqn:Ap} holds.  By
Lemma~\ref{lemma:rev-fac}, \eqref{eqn:Ap} holds if there
exist $\mu_1,\,\mu_2\in A_1$ such that
\[ 
H_1^{-\frac{\alpha p_0}{s}}H_2^{\frac{p_0}{q_0}}
w^{\frac{qp_0}{q_0}}\in A_{\frac{p_0}{p_-}}\cap RH_{\big(\frac{p_+}{p_0}\big)'}
=
W^{p_0} =
  \mu_2^{\frac{1}{(\frac{p_+}{p_0})'}}\mu_1^{1-\frac{p_0}{p_-}}. \]
By the $A_1$ property of the Rubio de Francia iteration algorithms, we
have that
\begin{gather*}
  \mu_1 = H_1^{\frac{q}{\tau}}
  w^{\frac{q}{\tau}-\frac{p}{\tau}\big(\frac{p_+}{p}\big)'} =
  \R_1(h_1^\delta
  w^\epsilon)\in A_1, \\
\mu_2 = H_2^{\frac{1}{\tau'}\big(\frac{q}{s}\big)'}
w^{\frac{\sigma}{\tau'}+\frac{q}{\tau'}} =\R_2(h^\beta
  w^\gamma)\in A_1.  
\end{gather*}
If we substitute these expressions into the above formula and equate
exponents, we see that equality holds if 
\begin{gather}
\label{eqn:exp1}
\frac{\alpha p_0}{s} = \frac{q}{\tau}\bigg(\frac{p_0}{p_-}-1\bigg), \\
\label{eqn:exp2}
\frac{p_0}{q_0}=\frac{1}{\tau'}\bigg(\frac{q}{s}\bigg)'\frac{1}{\big(\frac{p_+}{p_0}\big)'}, \\
\label{eqn:exp3}
\frac{qp_0}{q_0}
=
\bigg(\frac{\sigma}{\tau'}+\frac{q}{\tau'}\bigg)
\frac{1}{\big(\frac{p_+}{p_0}\big)'} +
\bigg(\frac{q}{\tau}-
\frac{p}{\tau}\bigg(\frac{p_+}{p}\bigg)'\bigg)
\bigg(1-\frac{p_0}{p_-}\bigg) 
. 
\end{gather}
If we use our choice of $\alpha$ on the left-hand side of
\eqref{eqn:exp1} and \eqref{eqn:s1} on the right-hand side, it is
straightforward to see that \eqref{eqn:exp1} holds. Additionally, if
we use \eqref{eqn:s2} on the right-hand side of \eqref{eqn:exp2}, we
see that the latter also holds. (It was the necessity of these
two identities for the proof that is the reason for our original
choice of $s$.)  To show that \eqref{eqn:exp3} holds, note that by
\eqref{eqn:tau-prime} and our choice of $\sigma$ we have that
\[ \frac{\sigma}{\tau'} 
= 
\frac{p}{\frac{p}{p_-}-1} \frac{\frac{1}{p_-}-\frac{1}{p}}{\frac{1}{p_-}-\frac{1}{p_+}}
= 
\frac{1}{\frac{1}{p_-}-\frac{1}{p_+}}. \]
Given this we can expand the right-hand side of \eqref{eqn:exp3}:
\begin{align*}
& \bigg(\frac{\sigma}{\tau'}+\frac{q}{\tau'}\bigg)
\frac{1}{\big(\frac{p_+}{p_0}\big)'} +
\bigg(\frac{q}{\tau}-\frac{p}{\tau}\bigg(\frac{p_+}{p}\bigg)'
\bigg)\bigg(1-\frac{p_0}{p_-}\bigg)
  \\
& \qquad = 
\bigg(\frac{1}{\frac{1}{p_-}-\frac{1}{p_+}}+q
  \frac{\frac{1}{p_-}-\frac{1}{p}}{\frac{1}{p_-}-\frac{1}{p_+}} \bigg)
p_0\bigg(\frac{1}{p_0}-\frac{1}{p_+}\bigg) \\
& \qquad \qquad \qquad +
\bigg(q \frac{\frac{1}{p}-\frac{1}{p_+}}{\frac{1}{p_-}-\frac{1}{p_+}}-
\frac{1}{\frac{1}{p_-}-\frac{1}{p_+}}\bigg)
p_0\bigg(\frac{1}{p_0}-\frac{1}{p_-}\bigg) \\
& \qquad =\frac{p_0}{\frac{1}{p_-}-\frac{1}{p_+}}
\bigg[ \frac{1}{p_-}-\frac{1}{p_+}
  +\frac{q}{p}\bigg(\frac{1}{p_+}-\frac{1}{p_-}\bigg) 
+\frac{q}{p_0}\bigg(\frac{1}{p_-}-\frac{1}{p_+}\bigg)\bigg] \\
& \qquad =  qp_0\bigg[ \frac{1}{q}-\frac{1}{p}+\frac{1}{p_0}\bigg] \\
& \qquad = \frac{qp_0}{q_0}.
\end{align*}
This completes the proof of Case I. \qed

\medskip

\subsection*{Proof of Theorem~\ref{thm:off-diag-limited}. Case II: $\boldsymbol{p_0=p_-}$}

Fix $p_-<p<p_+$ and $w$ such that
$w^p\in A_{\frac{p}{p_-}} \cap RH_{\big(\frac{p_+}{p}\big)'}$ and note
that in this case $p_-=p_0>0$ and $q_-=q_0$ by \eqref{eq:defqpm}. The
proof is similar to the proof of Case I and we
indicate the main changes. First, in
this case \eqref{eqn:s1} gives $s=q_0>0$. Thus, $s=q_0=q_-<q$ by \eqref{eq:defqpm} and the
fact that $p_-<p$.  Furthermore \eqref{eqn:s2} holds in this case. 

We now argue as before, but in this case we do not need to introduce
$H_1$. 
Since $s<q$, by rescaling and duality we have that
\[ 
\|f\|_{L^q(w^q)}^s = \|f^s\|_{L^{\frac{q}{s}}(w^q)} 
= 
\int_\rn f^s h_2 w^q\,dx
\le
\int_\rn f^s H_2 w^q\,dx
, \]
where $h_2$ is a non-negative function in $L^{(\frac{q}{s})'}(w^q)$
with $\|h_2\|_{L^{(\frac{q}{s})'}(w^q)}=1$ and $H_2$ is such that
$h_2\leq H_2$; we will determine the exact value below.  If we assume
further that
$\|H_2\|_{L^{\big(\frac{q}{s}\big)'}(w^q)}\leq C_2<\infty$, it follows
by assumption that 
$$
\int_\rn f^s H_2 w^q\,dx\le
\|f^s\|_{L^{\frac{q}{s}}(w^q)} \|H_2\|_{L^{\big(\frac{q}{s}\big)'}(w^q)}
\le
C_2\|f\|_{L^q(w^q)}^s<\infty.
$$
Define
$\varphi=\big(\frac{q}{s}\big)'\frac{q_0}{p_0}=\big(\frac{q}{q_0}\big)'\frac{q_0}{p_0}$;
then we have that
\[ \frac{1}{\varphi} =
  \frac{p_0}{q_0}\frac{1}{\big(\frac{q_0}{s}\big)'}=
  \frac{p_0}{q_0}\bigg(1-\frac{q_0}{q}\bigg)=
  p_0\bigg(\frac{1}{q_0}-\frac{1}{q}\bigg) =
  p_0\bigg(\frac{1}{p_0}-\frac{1}{p}\bigg) =1-\frac{p_0}{p} < 1, \]
which implies that $\varphi'=\frac{p}{p_0}$.  Now let
$W^{q_0}=H_2 w^q$ and assume that
$W^{p_0}\in A_{\frac{p_0}{p_-}} \cap
RH_{\big(\frac{p_+}{p_0}\big)'}=A_1\cap
RH_{\big(\frac{p_+}{p_-}\big)'}$, or equivalently (by Lemma
\ref{lemma:cjn}), $W^{p_0\big(\frac{p_+}{p_-}\big)'}\in A_1$. Then by 
our hypothesis \eqref{eqn:off-diag-limited1} we get
\begin{align*}
\|f\|_{L^q(w^q)}^s
& = \int_\rn f^{q_0} W^{q_0}\,dx  \\
& \leq C\left(\int_\rn g^{p_0} W^{p_0}\,dx\right)^{\frac{q_0}{p_0}} \\
& =  C\left(\int_\rn g^{p_0} 
H_2^{\frac{p_0}{q_0}}w^{q\frac{p_0}{q_0}}w^{-q}w^q\,dx\right)^{\frac{q_0}{p_0}}
  \\ 
& 
\leq 
\left(\int_\rn g^{p_0\varphi'} 
w^{q(\frac{p_0}{q_0}-1)\varphi'}w^q\,dx\right)^{\frac{q_0}{\varphi' p_0}}
\left(\int_\rn  H_2^{\big(\frac{q}{s}\big)'}w^q\,dx\right)^{\frac{q_0}{\varphi p_0} }.
\\
&\le
C_2
\left(\int_\rn g^{p_0\varphi'} 
w^{q(\frac{p_0}{q_0}-1)\varphi'}w^q\,dx\right)^{\frac{q_0}{\varphi' p_0}}
\\
&=
C_2
\left(\int_\rn g^p 
w^p\,dx\right)^{\frac{q_0}{p_0}},
\end{align*}
where in the last equality we have used that
$$
q\Big(\frac{p_0}{q_0}-1\Big)\varphi'+q
=
q\Big(\frac{p_0}{q_0}-1\Big)\frac{p}{p_0}+q
=
qp\Big(\frac1{q_0}-\frac1{p_0}\Big)+q
=
qp\Big(\frac1{q}-\frac1{p}\Big)+q
=p.
$$

\medskip

Therefore, to complete the proof we need to show that we can construct
a non-negative function $H_2$ such that 
\begin{gather}
\label{eqn:H2-norm-case2}
\|H_2\|_{L^{(\frac{q}{s})'}(w^q)} \leq C_2, \\
\label{eqn:H2-pt-case2}
h_2 \leq H_2;
\end{gather}
and such that  the weight 
$W= H_2^{\frac{1}{q_0}}  w^{\frac{q}{q_0}}$
satisfies
\begin{equation} \label{eqn:Ap-case2}
W^{p_0\big(\frac{p_+}{p_-}\big)'}
= 
H_2^{\frac{p_0}{q_0}\big(\frac{p_+}{p_-}\big)'}
w^{\frac{qp_0}{q_0}\big(\frac{p_+}{p_-}\big)'}
\in A_1. 
\end{equation}

We construct $H_2$ exactly as in the proof of Case I, and as before we
have \eqref{eqn:H2-norm-case2} and \eqref{eqn:H2-pt-case2}. It remains
to show \eqref{eqn:Ap-case2}. By \eqref{eqn:s2},
$$
\frac{1}{\beta}\frac{p_0}{q_0}\Big(\frac{p_+}{p_-}\Big)'
=
\frac{\tau'}{\big(\frac{q}{s}\big)'}\frac{p_0}{q_0}\Big(\frac{p_+}{p_0}\Big)'
=
\frac{1}{\big(\frac{q}{s}\big)'}\frac{q-s}{s}=1.
$$ 
On the other hand, recalling that $p_0=p_-$ and $s=q_0$ we obtain
$$
q-\frac{\gamma}{\beta}
=
q-\frac{\sigma+q}{\big(\frac{q}{s}\big)'}
=
q_0-\frac{p}{\big(\frac{q}{s}\big)'\big(\frac{p}{p_-}-1\big)}
=
q_0-q_0\frac{\frac{1}{q_0}-\frac1{q}}{\frac{1}{p_0}-\frac1{p}}
=
0.
$$
Thus,
\begin{multline*}
W^{p_0\big(\frac{p_+}{p_-}\big)'}
=
H_2^{\frac{p_0}{q_0}\big(\frac{p_+}{p_-}\big)'}
w^{\frac{qp_0}{q_0}\big(\frac{p_+}{p_-}\big)'} \\
=
\R_2(h_2^\beta  w^\gamma)^{\frac{1}{\beta}\frac{p_0}{q_0}\big(\frac{p_+}{p_-}\big)'}
w^{\frac{p_0}{q_0}\big(\frac{p_+}{p_-}\big)'\big(q-\frac{\gamma}{\beta}\big)}
=
\R_2(h_2^\beta  w^\gamma)\in A_1,
\end{multline*}
which concludes the proof of Case II. \qed

\medskip

\subsection*{Proof of Theorem~\ref{thm:off-diag-limited}. Case III: 
$\boldsymbol{p_0=p_+}$ and $\boldsymbol{p_->0}$}
Fix $p_-<p<p_+$ and $w$ such that
$w^p\in A_{\frac{p}{p_-}} \cap RH_{\big(\frac{p_+}{p}\big)'}$ and note
that in this case $p_+=p_0<\infty$ and $q_+=q_0$ by
\eqref{eq:defqpm}. We again follow the proof of Case I and we indicate
the main changes. First, if we
define $s$ as in \eqref{eqn:s1} and since \eqref{eqn:s2} is also valid
in this context,  then $0<s=q<q_+=q_0$ by \eqref{eq:defqpm}
and the fact that $p<p_+$. 

We now argue as before, but in this case we do not need to use duality
or introduce
$H_2$.  Since $s=q$, if we fix
$\alpha=\frac{s}{\big(\frac{q_0}{s}\big)'}$, then by H\"older's
inequality, 
\begin{align*}
\|f\|_{L^q(w^q)}^{q_0} 
&=
\left(\int_\rn f^s H_1^{-\alpha} H_1^{\alpha}w^q\,dx\right)^{\frac{q_0}{q}}
\\
&\leq 
\left(\int_\rn f^{q_0}H_1^{-\alpha\frac{q_0}{s}}w^q\,dx\right)
\left( \int_\rn H_1^{\alpha\left(\frac{q_0}{s}\right)'} w^q\,dx\right)^{\frac{q_0}{q\big(\frac{q_0}{s}\big)'}} 
\\
&\le
C_1^{{q_0}/{\big(\frac{q_0}{s}\big)'}} \int_\rn f^{q_0}H_1^{-\alpha\frac{q_0}{s}}w^q\,dx
, \end{align*}
where $0<H_1<\infty$ is in $L^q(w^q)$ with  $\|H_1\|_{L^{q}(w^q)}\le
C_1<\infty$. 
We will determine the exact value  below. If we also assume that 
$f\leq H_1\|f\|_{L^q(w^q)}$, then 
\begin{multline*}
\int_\rn f^{q_0}H_1^{-\alpha\frac{q_0}{s}}w^q\,dx
\leq \|f\|_{L^q(w^q)}^{q_0}\int_\rn H_1^{q_0} H_1^{-\alpha\frac{q_0}{s}}
  w^q\,dx \\
= \|f\|_{L^q(w^q)}^{q_0}\int_\rn H_1^s 
  w^q\,dx \le C_1^q \|f\|_{L^q(w^q)}^{q_0}< \infty. 
\end{multline*}
Thus, we can  apply~\eqref{eqn:off-diag-limited1} if we  let
$W^{q_0}= H_1^{-\alpha\frac{q_0}{s}}w^q$ and assume that
$W^{p_0}\in A_{\frac{p_0}{p_-}} \cap
RH_{\big(\frac{p_+}{p_0}\big)'}=A_{\frac{p_+}{p_-}} \cap RH_{\infty}$:
\begin{align*}
\|f\|_{L^q(w^q)}^{q_0} 
& 
\le
C_1^{{q_0}/{\big(\frac{q_0}{s}\big)'}}
 \int_\rn f^{q_0} W^{q_0}\,dx  \\
& 
\leq 
C\left(\int_\rn g^{p_0} W^{p_0}\,dx\right)^{\frac{q_0}{p_0}} \\
& 
=  
C\left(\int_\rn g^{p_0} H_1^{-\alpha \frac{p_0}{s}}
w^{q\frac{p_0}{q_0}}w^{-q}w^q\,dx\right)^{\frac{q_0}{p_0}}
  \\ 
& \leq 
C
\|g\|_{L^p(w^p)}^{q_0}
\left(\int_\rn H_1^q w^q\,dx\right)^{\frac{q_0}{p_0}}
\le
C C_1^{\frac{q_0q}{p_0}}
\|g\|_{L^p(w^p)}^{q_0},
\end{align*}
provided $H_1$ satisfies
\[ g^{p_0} H_1^{-\alpha \frac{p_0}{s}}
w^{q(\frac{p_0}{q_0}-1)} \leq H_1^q
\|g\|_{L^p(w^p)}^{p_0}. \]

\medskip

To complete the proof we need to show that we can construct
$H_1$ such that
\begin{gather}
\label{eqn:H1-norm-case3}
\|H_1\|_{L^q(w^q)} \leq C_1, \\
\label{eqn:H1-pt-case3}
g^{p_0} H_1^{-\alpha \frac{p_0}{s}}
w^{q(\frac{p_0}{q_0}-1)} \leq H_1^q \|g\|_{L^p(w^p)}^{p_0}, \\
\label{eqn:H1-f-case3}
0<H_1<\infty, \quad f\leq H_1\|f\|_{L^q(w^q)},
\end{gather}
and such that  the weight 
$W= H_1^{-\frac{\alpha}{s}}  w^{\frac{q}{q_0}}$
satisfies
\begin{equation} \label{eqn:Ap-case3}
W^{p_0}=
H_1^{-\frac{\alpha p_0}{s}}
w^{\frac{qp_0}{q_0}}\in A_{\frac{p_0}{p_-}}\cap RH_{\infty}. 
\end{equation}
Since $\alpha\frac{p_0}{s}=\frac{p_0}{q_0}(q_0-s)$, 
\eqref{eqn:H1-pt-case3} is equivalent to
\begin{equation} \label{eqn:H1-pt2-case3}
 g^{p_0}w^{q(\frac{p_0}{q_0}-1)} 
\leq
H_1^{q+\frac{p_0}{q_0}(q_0-s)}\|g\|_{L^p(w^p)}^{p_0}. 
\end{equation}
Using the fact that
$\frac{1}{p}-\frac{1}{q}=\frac{1}{p_0}-\frac{1}{q_0}$, and that $s=q$ we have that
$$
q +\frac{p_0}{q_0}(q_0-s)
=
q+p_0q\bigg(\frac{1}{q}-\frac{1}{q_0}\bigg)
= 
q+p_0q\bigg(\frac{1}{p}-\frac{1}{p_0}\bigg)
=q\frac{p_0}{p}.
$$
Similarly, we have that 
\[ q\bigg(\frac{p_0}{q_0}-1\bigg)
= qp_0\bigg(\frac{1}{q_0}-\frac{1}{p_0}\bigg)
=qp_0\bigg(\frac{1}{q}-\frac{1}{p}\bigg)
=
p_0\bigg(1-\frac{q}{p}\bigg)
. \]
Therefore, \eqref{eqn:H1-pt2-case3} (and hence \eqref{eqn:H1-pt-case3}) is equivalent to 
\begin{equation} \label{eqn:H1-pt3-case3}
 g^{\frac{p}{q}} w^{\frac{p}{q}-1} \leq H_1 \|g\|_{L^p(w^p)}^{\frac{p}{q}}. 
\end{equation}
We now construct $H_1$ exactly as in the proof of Case I, and we obtain
as before \eqref{eqn:H1-f-case3}, \eqref{eqn:H1-pt3-case3}, and
\eqref{eqn:H1-norm-case3}. It remains to show
\eqref{eqn:Ap-case3}. By \eqref{eqn:s1}
$$
\frac{\alpha p_0}{\delta s}
=
\frac{\tau p_0}{q \big(\frac{q_0}{s}\big)'}
=
\frac{1}{\big(\frac{q_0}{s}\big)'}
\frac{\frac{p_0}{p_-}-1}{\frac{q_0-s}{q_0}}
=
\frac{p_0}{p_-}-1,
$$
and also, since $p_0=p_+$,
\begin{multline*}
\frac{\epsilon\alpha p_0}{\delta s}
=
\Big(\frac{p_0}{p_-}-1\Big)\epsilon
=
\Big(\frac{p_0}{p_-}-1\Big)q\bigg(\frac{1}{\tau}-\frac{p\big(\frac{p_+}{p})'}{\tau}\bigg)
\\
=
p_0\Big(\frac1{p_-}-\frac1{p_0}\Big)\bigg(
q\frac{\frac{1}{p}-\frac{1}{p_+}}{\frac{1}{p_-}-\frac{1}{p_+}}
-
\frac{1}{\frac{1}{p_-}-\frac{1}{p_+}}\bigg)
=
qp_0\Big( \frac1{p}-\frac1{p_0}-\frac1q\Big)
=
-\frac{qp_0}{q_0}.
\end{multline*}
Together, these  imply that 
$$
W^{p_0}=
H_1^{-\frac{\alpha p_0}{s}}
w^{\frac{qp_0}{q_0}}
=
 \R_1(h_1^\delta
  w^\epsilon)^{1-\frac{p_0}{p_-}}
	\in 
	A_{\frac{p_0}{p_-}}\cap RH_{\infty};
$$
the inclusion follows from  Lemma \ref{lemma:rev-fac} and the fact that $\R_1(h_1^\delta
  w^\epsilon)\in A_1$. This completes the proof of Case III. \qed

\medskip

\subsection*{Proof of Theorem~\ref{thm:off-diag-limited}. 
Case IV: $\boldsymbol{p_-=0}$ and $\boldsymbol{p_-<p_0\le p_+}$}

In this case we adapt ideas from \cite[Section
3.1]{martell-prisuelos}. Fix $p$, $q$ such that $0=p_-<p<p_+$,
$0<q<\infty$ and
$\frac{1}{p}-\frac{1}{q}=\frac{1}{p_0}-\frac{1}{q_0}$, and let $v$ be
such that
$v^{p}\in A_{\frac{p}{p_-}} \cap RH_{\big(\frac{p_+}{p}\big)'}=
RH_{\big(\frac{p_+}{p}\big)'}$. Since $RH_1=A_\infty$, there exists
$0<\epsilon<\min\{p_0,p\}$ such that
$v^p\in A_{\frac{p}{\epsilon}}$. Set $\widetilde{p}_-=\epsilon>0$; then
$\widetilde{p}_-<p_0\le p_+$ and 
\eqref{eqn:off-diag-limited1} holds for all
$w^{p_0}\in A_{\frac{p_0}{\widetilde{p}_-}} \cap
RH_{\big(\frac{p_+}{p_0}\big)'} \subset
RH_{\big(\frac{p_+}{p_0}\big)'}= A_{\frac{p_0}{p_-}} \cap
RH_{\big(\frac{p_+}{p_0}\big)'}$.  
Thus, we can use Cases I and III
with $\widetilde{p}_->0$ in place of $p_-$ to  conclude that
\eqref{eqn:off-diag-limited2} holds for every $\widetilde{p}$,
$\widetilde{q}$ such that $\widetilde{p}_-<\widetilde{p}<p_+$,
$0<\widetilde{q}<\infty$ and
$\frac{1}{\widetilde{p}}-\frac{1}{\widetilde{q}}=\frac{1}{p_0}-\frac{1}{q_0}$,
and every weight $w$ such that
$w^{\widetilde{p}}\in A_{\frac{\widetilde{p}}{\widetilde{p}_-}} \cap
RH_{\big(\frac{p_+}{\widetilde{p}}\big)'}$. 
If we take $\widetilde{p}=p$,
$\widetilde{q}=q$ and $w=v$, our choice of $\epsilon$ guarantees that
$\widetilde{p}_-=\epsilon<p<p_+$, $0<q<\infty$ and
$\frac{1}{p}-\frac{1}{q}=\frac{1}{p_0}-\frac{1}{q_0}$. Moreover,
$v^p\in A_{\frac{p}{\epsilon}}\cap RH_{\big(\frac{p_+}{p}\big)'}=
A_{\frac{\widetilde{p}}{\widetilde{p}_-}} \cap
RH_{\big(\frac{p_+}{\widetilde{p}}\big)'}$. Thus,
\eqref{eqn:off-diag-limited2} holds and the proof of Case IV is complete. \qed

\medskip

\subsection*{Proof of Theorem \ref{thm:extrapol-multi}}
Our proof of Theorem \ref{thm:extrapol-multi} is a modification of the
proof of  multilinear extrapolation in~\cite[Theorem~6.1]{MR2754896}.  We include
the details so that we can explain the use of families of
extrapolation pairs.  The essential idea is to reduce the problem to a
linear one by acting on one function at a time.

For $2\leq j \leq m$, fix weights
$w_j$ such that
$w_j^{p_j} \in A_{\frac{p_j}{r_j^-}}\cap
RH_{\big(\frac{r_j^+}{p_j}\big)'}$.   Fix functions $f_j$, $2\leq j \leq m$, such that there exists
functions $f$ and $g$ with $(f,g,f_2,\ldots,f_m)\in \F$.  Assume that for
each $j$, $0<\|f_j\|_{L^{p_j}(w_j^{p_j})}<\infty$.   (We will remove
this restriction below.)   Define the new family of extrapolation
pairs
\[ \F_1 = 
\{ (F,g)= (f\prod_{j=2}^m w_j\|f_j\|_{L^{p_j}(w_j^{p_j})}^{-1},  g)
 : (f,g, f_2,\ldots,f_m)\in \F \}.  \]
If $f \in L^p(w^p)$, then $F\in L^{p}(w_1^{p})$, so by our hypothesis~\eqref{eq:extrap-multi-hyp},
\begin{equation}\label{eq:multi-1-coor}
\|F\|_{L^{p}(w_1^p)} \leq C\|g\|_{L^{p_1}(w_1^{p_1})}
\end{equation}
for all
$w_1^{p_1} \in A_{\frac{p_1}{r_1^-}}\cap
RH_{\big(\frac{r_1^+}{p_1}\big)'}$.  Note that $p<p_1$ and so
$\frac1p-\frac1{p_1}+\frac1{r_1^+}>0$.
Therefore, by Theorem~\ref{thm:off-diag-limited}, for all pairs
$(F,g)\in \F_1$ with  $\|F\|_{L^q(w_1^q)}<\infty$, and 
for all $r_1^- <q_1 < r_1^+$ and
all $w_1^{q_1} \in  A_{\frac{q_1}{r_1^-}}\cap
RH_{\big(\frac{r_1^+}{q_1}\big)'}$, 
\[ \|F\|_{L^q(w_1^q)} \leq 
C\| g\|_{L^{q_1}(w_1^{q_1})}, \]
where $\frac{1}{q}-\frac{1}{q_1}=\frac{1}{p}-\frac{1}{p_1}$ and so $
\frac{1}{q} = \frac{1}{q_1}+ \sum_{j=2}^m \frac{1}{p_j}$.
Therefore, by our definition of $F$, $\|f\|_{L^q(w^q)} <\infty$ we can
rewrite this as
\[
\|f\|_{L^q(w^q)} 
\leq C\| g\|_{L^{q_1}(w_1^{q_1}) }
\prod_{j=2}^m \|f_j\|_{L^{p_j}(w_1^{p_j})}. 
\]
This inequality still holds even if we  remove the restriction $0<\|f_j\|_{
    L^{p_j}(w_j^{p_j})} <\infty$.   If for some $j$, $\|f_j\|_{
    L^{p_j}(w_j^{p_j})}=\infty$, this inequality clearly holds; if  $\|f_j\|_{
    L^{p_j}(w_j^{p_j})}=0$, then~\eqref{eq:extrap-multi-hyp} implies
  that $f=F=0$, and this inequality again holds.

We can repeat this argument for any such collection of $f_j$, $2\leq j
\leq m$.  Therefore, we have shown that for all $(f_1,\ldots,f_m)\in
\F$ with $f\in L^q(w^q)$,
\[
\|f\|_{L^q(w^q)} 
\leq C\| f_1\|_{L^{q_1}(w_1^{q_1}) }
\prod_{j=2}^m \|f_j\|_{L^{p_j}(w_1^{p_j})}.
\]

To complete the proof, fix
$f_1,f_3,\ldots f_m$, and repeat the above argument in the second
coordinate, etc.  Then by induction we get the desired conclusion.

\medskip

We now prove the vector-valued inequalities
\eqref{eq:extrap-multi-con-vv}. The extension of scalar inequalities
to vector-valued inequalities via extrapolation is well-known in the
linear case: see~\cite[Corollary~3.12]{cruz-martell-perezBook}.  The
argument is nearly the same in the multilinear setting.  Fix $s_j$,
$r_j^-<s_j<r_j^+$, for $1\le j\le m$ and set
$\displaystyle \frac{1}{s} = \sum_{j=1}^m \frac{1}{s_j}$. Define a new
family
\begin{multline*}
 \widetilde{\F} 
= 
\bigg\{ (F,F_1,\ldots, F_m)
= \bigg( \bigg(\sum_k (f^k)^s\bigg)^{\frac{1}{s}}, 
\bigg(\sum_k (f_1^k)^{s_1}\bigg)^{\frac{1}{s_1}}, \ldots,
\bigg(\sum_k (f_m^k)^{s_m}\bigg)^{\frac{1}{s_m}} \bigg) : 
\\
\big\{(f^k,f_1^k,\ldots, f_m^k)\}_k\subset \F\bigg\}.  
\end{multline*}
Without loss of generality we may  assume that all of the
sums in the definition of $\widetilde{\F}$ are finite; the conclusion
for infinite sums follows by the monotone convergence theorem.
Then, given any collection of weights
$w_1,\ldots,w_m$ with $w_j^{s_j} \in A_{\frac{s_j}{r_j^-}}\cap
RH_{\big(\frac{r_j^+}{s_j}\big)'}$ and $w=w_1\cdots w_m$, 
if $\|F\|_{L^s(w^s)}<\infty$, then by \eqref{eq:extrap-multi-con} we
have that
\begin{multline}\label{eq:extrap-multi-con-vv:proof}
\|F\|_{L^s(w^s)}
=
\bigg( \sum_k \|f^k\|_{L^s(w^s)}^s\bigg)^{\frac1s} 
\leq
C
\bigg(\sum_k\prod_{j=1}^m 
\|f_j^k\|_{L^{s_j}(w_j^{s_j})}^s
\bigg)^{\frac1s}
\\
\leq
C
\prod_{j=1}^m \bigg(\sum_k 
\|f_j^k\|_{L^{s_j}(w_j^{s_j})}^{s_j}
\bigg)^{\frac1{s_j}}
=
C
\prod_{j=1}^m \|F_j\|_{L^{s_j}(w_j^{s_j})},
\end{multline}
where in the second estimate we used H\"older's inequality with
respect to sums. We can now apply the first part of
Theorem~\ref{thm:extrapol-multi} to $\widetilde{F}$,  where we use
\eqref{eq:extrap-multi-con-vv:proof} for the initial
estimate in place of \eqref{eq:extrap-multi-hyp}.  We thus get
\begin{equation} \label{eqn:vv-est}
\|F\|_{L^q(w^q)}
\le
C\prod_{j=1}^m \|F_j\|_{L^{q_j}(w_j^{q_j})}.
\end{equation}
for all exponents $q_j$, $r_j^-<q_j<r_j^+$, all weights
$w_j^{q_j} \in A_{\frac{q_j}{r_j^-}}\cap
RH_{\big(\frac{r_j^+}{q_j}\big)'}$, $w=w_1\cdots w_m$, and
$\displaystyle \frac{1}{q} = \sum_{j=1}^m
\frac{1}{q_j}$. Inequality~\eqref{eqn:vv-est} holds for all
$(F,F_1,\ldots,F_m)\in\widetilde{\F}$ for which
$\|F\|_{L^q(w^q)}<\infty$. But this is exactly
\eqref{eq:extrap-multi-con-vv} and the proof is complete. \qed

\medskip

\subsection*{Proof of Corollaries \ref{cor:extrapol-mult} and \ref{corol:extrapol-rest}}
We will prove
Corollary~\ref{corol:extrapol-rest}; the proof of
Corollary~\ref{cor:extrapol-mult} is identical. 
The proof follows as in \cite[Section
3.1]{martell-prisuelos}.  
Given a family of extrapolation pairs $\mathcal{F}$ as in the
statement and any $N>0$, define  the new
family 
$$
\mathcal{F}_N
:=
\big\{(f_N,g): (f,g)\in\mathcal{F}, f_N:=f\chi_{\{x\in B(0,N):f(x)\leq N\}}
\big\}.
$$
Note that for all $0<r<\infty$ and $w^r\in A_\infty$, 
\begin{equation}\label{eq:LHS-fN}
\|f_N\|_{L^r(w^r)}^r\leq  N^r w^r(B(0,N))<\infty.
\end{equation}
Since $f_N\le f$, by our hypothesis we get that
\eqref{eqn:off-diag-limited1} holds for every pair in $\mathcal{F}_N$
(with a constant independent of $N$) with a left-hand side that is
always finite by~\eqref{eq:LHS-fN} and
Remark~\ref{remark:rest-expon:0}.  Therefore, we can apply
Theorem~\ref{thm:off-diag-limited} to $\mathcal{F}_N$ to conclude that
\eqref{eqn:off-diag-limited2} holds for every pair
$(f_N,g)\in\mathcal{F}_N$ (with a constant that is again independent
of $N$), since again the left-hand side is always finite. The desired
inequality follows at once if we let $N\to\infty$ and apply the
monotone convergence theorem. \qed


\medskip

\section{Proofs of the applications}
\label{section:BHT}

We now prove Theorems~\ref{thm:bht-DCU-JMM}, \ref{thm:bht-vector}, and
\ref{theorem:CZ-MZ}, and 
Corollary~\ref{cor:bht-A1}.  We also sketch the ideas needed to prove
the result in Remark~\ref{remark:iteration}.

\subsection*{Proof of Theorem~\ref{thm:bht-DCU-JMM}}
We start with the first part of the theorem. Let
$p_1, p_2\in (1,\infty)$ be such that
$\frac{1}{p}=\frac{1}{p_1}+\frac{1}{p_2}<1$,  fix
$w_1^{2p_1}\in A_{p_1}$, $w_2^{2p_2}\in A_{p_2}$, and let $w=w_1w_2$.
Then by Theorem~\ref{thm:bht-cdo}, $BH : L^{p_1}(w_1^{p_1}) \times
L^{p_2}(w_2^{p_2}) \rightarrow L^p(w^p)$.  By Lemma~\ref{lemma:cjn}, 
$w_i^{2p_i}\in A_{p_i}$ if and only if
$w^{p_i}\in A_{\frac{p_i+1}{2}}\cap RH_2$. Thus, if we set
$r_i^-=\frac{2p_i}{p_i+1}$ and $r_i^+=2p_i$, then
$1<r_i^-<p_i<r_i^+<\infty$ and
$w_i^{p_i}\in A_{\frac{p_i}{r_i^-}}\cap
RH_{\big(\frac{r_i^+}{p_i}\big)'}$. We can then apply
Corollary~\ref{cor:extrapol-mult} to the family 
\[  \F = \{ (|B(f,g)|,|f|,|g|) : f, g \in L^\infty_c \} \]
to conclude that for
all $r_i^-<q_i<r_i^+$ and
$w_i^{q_i}\in A_{\frac{q_i}{r_i^-}}\cap
RH_{\big(\frac{r_i^+}{q_i}\big)'}$, the bilinear Hilbert transform
$BH$ is bounded from $L^{q_1}(w_1^{q_1}) \times L^{q_2}(w_2^{q_2})$ into
$L^q(w^q)$ where $\frac1q=\frac{1}{q_1}+\frac{1}{q_2}$ and
$w=w_1w_2$.  (Here we use the fact that $L^\infty_c$ is dense any
space $L^r(w^r)$ if $w^r$ is locally integrable, and the fact that
$BH$ is bilinear to extend the inequality on triples in $\F$ that we get from
Theorem~\ref{thm:extrapol-multi} to all of $L^{q_1}(w_1^{q_1}) \times
L^{q_2}(w_2^{q_2})$.) 
Again by Lemma~\ref{lemma:cjn}, 
the conditions on the weights are equivalent to
$w_i^{2r_i}\in A_{r_i}$, where
$r_i=\big(\frac2{q_i}-\frac1{p_i}\big)^{-1}$. Note that $1<r_i<\infty$
since $r_i^-<q_i<r_i^+$. This completes the proof of the first part of
Theorem~\ref{thm:bht-DCU-JMM}.

\medskip

To prove the second part of the theorem, fix $1<q_1,q_2<\infty$ such that
$\frac{1}{q}=\frac{1}{q_1}+\frac{1}{q_2}<\frac32$. We want to use
the previous argument:  therefore,  we need to find
$p_1, p_2\in (1,\infty)$ such that
$\frac{1}{p}=\frac{1}{p_1}+\frac{1}{p_2}<1$ and $r_i^-<q_i<r_i^+$, where
\begin{equation}\label{eq:formula-r_ipm}
\frac1{r_i^+}=\frac1{2p_i}
<\frac1{q_i}
<
\frac1{r_i^-}=\frac1{2p_i}+\frac12.
\end{equation}
Since $1<p_1,p_2<\infty$, this can be rewritten as 
\begin{equation}\label{eq:formula-r_ipm:equiv}
0\le 2\left(\max\left\{\frac12,\frac1{q_i}\right\}-\frac12\right)
<
\frac1{p_i}
<
2\min\left\{\frac12,\frac1{q_i}\right\}\le 1.
\end{equation}
Before choosing $p_1,p_2$ we claim that
\begin{equation}\label{eq:restriction}
\sum_{i=1}^2 \max\left\{\frac12,\frac1{q_i}\right\}<\frac32.
\end{equation}
To see that this holds, note that 
$$
\sum_{i=1}^2 \max\left\{\frac12,\frac1{q_i}\right\}
=
\left\{
\begin{array}{ll}
1&\text{if } \max\{\frac1{q_1},\frac1{q_2}\}\le\frac12,
\\[8pt]
\frac12+\max\{\frac1{q_1},\frac1{q_2}\}
& \text{if }	
\min\{\frac1{q_1},\frac1q_2\}\le\frac12
\le\max\{\frac1{q_1},\frac1{q_2}\},
\\[8pt]
\frac1{q_1}+\frac1{q_2}
&
\text{if }	
\min\{\frac1{q_1},\frac1{q_2}\}>\frac12.
\end{array}
\right.
$$
and in every case this is strictly smaller than $\frac32$ since $q_1,q_2>1$ and $\frac{1}{q}=\frac{1}{q_1}+\frac{1}{q_2}<\frac32$.

Now define
\begin{equation}\label{eq:choice-pi}
\frac1{p_i}:=2\left(\max\left\{\frac12,\frac1{q_i}\right\}-\frac12+\eta_i\right),
\qquad i=1,2,
\end{equation}
where we fix $\eta_1,\,\eta_2>0$ so that 
\begin{equation}\label{eq:conds-eta}
\eta_1+\eta_2<
\frac32-\sum_{i=1}^2 \max\left\{\frac12,\frac1{q_i}\right\}
\quad\text{and}\quad 
0<\eta_i<\min\left\{\frac1{q_i},\frac1{q_i'}\right\}, \quad i=1,2.
\end{equation}
That we can find such $\eta_1,\eta_2$ follows
from~\eqref{eq:restriction}. (As  will be clear from the proof, we can choose $\eta_i$
as close to $0$ as we want; we will use this fact in the proof of
Corollary~\ref{cor:bht-A1} below.)

With this choice we claim that 
\eqref{eq:formula-r_ipm:equiv} holds and also that
$\frac{1}{p}=\frac{1}{p_1}+\frac{1}{p_2}<1$. We first prove the latter
inequality:  by the first condition in \eqref{eq:conds-eta},
$$
\frac{1}{p}
=
\frac{1}{p_1}+\frac{1}{p_2}
=
-2+2\,\sum_{i=1}^2
\max\left\{\frac12,\frac1{q_i}\right\}+2\sum_{i=1}^2 \eta_i< 1.
$$
To prove \eqref{eq:formula-r_ipm:equiv} we first observe that since $\eta_i>0$,
$$
2\left(\max\left\{\frac12,\frac1{q_i}\right\}-\frac12\right)
<
2\left(\max\left\{\frac12,\frac1{q_i}\right\}-\frac12+\eta_i\right)
=
\frac1{p_i}.
$$
To obtain the other half of \eqref{eq:formula-r_ipm:equiv} we consider two cases. If $\max\{\frac12,\frac1{q_i}\}=\frac12$, then
$$
\frac1{p_i}
=
2\eta_i
<
\frac2{q_i}
=
2\min\left\{
\frac12,
\frac1{q_i}\right\}.
$$
On the other hand, if $\max\{\frac12,\frac1{q_i}\}=\frac1{q_i}$, then
$$
\frac1{p_i}=\frac2{q_i}-1+2\eta_i
<
\frac2{q_i}-1+\frac2{q_i'}
=
1
=
2\min\left\{
\frac12,\frac1{q_i}\right\}.
$$
This completes the proof of \eqref{eq:formula-r_ipm:equiv} and hence
the proof of Theorem~\ref{thm:bht-DCU-JMM}.\qed 

\medskip

\subsection*{Proof of Corollary~\ref{cor:bht-A1}}
This result follows by considering more carefully the proof of
Theorem~\ref{thm:bht-DCU-JMM}.  Fix
$1<q_1,q_2<\infty$ such that
$\frac{1}{q}=\frac{1}{q_1}+\frac{1}{q_2}<\frac32$ and
$w_i^{q_i}\in A_{\max\{1,\frac{q_i}2\}}\cap
RH_{\max\{1,\frac2{q_i}\}}$.  We now choose $p_i$ as in
\eqref{eq:choice-pi} and \eqref{eq:conds-eta}, though below  we will take
$\eta_i$ much smaller. As we showed above, $1<p_1,p_2<\infty$,
$\frac{1}{p}=\frac{1}{p_1}+\frac{1}{p_2}<1$, and
\eqref{eq:formula-r_ipm:equiv} holds. Hence,
\eqref{eq:formula-r_ipm} holds and so by the first part of Theorem
~\ref{thm:bht-DCU-JMM}, we get that the bilinear Hilbert transform is bounded from
$L^{q_1}(u_1^{q_1}) \times L^{q_2}(u_2^{q_2})$ into $L^q(u^q)$ where
$\frac1q=\frac{1}{q_1}+\frac{1}{q_2}$ and $u=u_1u_2$, for all
$u_i^{q_i}\in A_{\frac{q_i}{r_i^-}}\cap
RH_{\big(\frac{r_i^+}{q_i}\big)'}$ with
$$
\frac{q_i}{r_i^-}
=
q_i\left(\frac1{2p_i}+\frac12\right)
=
q_i\left(\max\left\{\frac12,\frac1{q_i}\right\}+\eta_i\right)
=
\max\left\{1,\frac{q_i}2\right\}+q_i\eta_i,
$$
and
$$
\frac1{\big(\frac{r_i^+}{q_i}\big)'}
=
1-\frac{q_i}{r_i^+}
=
1-\frac{q_i}{2p_i}
=
1-q_i\left(\max\left\{\frac12,\frac1{q_i}\right\}-\frac12+\eta_i\right)
=
\min\left\{1,\frac{q_i}2\right\}-q_i\eta_i.
$$
Note that $w_i^{q_i}\in A_{\max\{1,\frac{q_i}2\}}$ immediately implies that
$w_i^{q_i}\in A_{\frac{q_i}{r_i^-}}$. On the other hand, since
$w_i^{q_i}\in RH_{\max\{1,\frac2{q_i}\}}$, by the openness of the
reverse H\"older classes we can find $0<\theta<1$ close to $1$ such that
$w_i^{q_i}\in RH_{\frac1{\theta}\max\{1,\frac2{q_i}\}}$.  Therefore,
in choosing the $\eta_i$ we
assume that \eqref{eq:conds-eta} holds and  that
$0<\eta_i<(1-\theta)\min\left\{1,\frac{q_i}2\right\}$.  But then
$$
\frac1{\big(\frac{r_i^+}{q_i}\big)'}
=
\min\left\{1,\frac{q_i}2\right\}-q_i\eta_i
>
\min\left\{1,\frac{q_i}2\right\}-q_i(1-\theta)\min\left\{\frac12,\frac1{q_i}\right\}
=
\theta\min\left\{1,\frac{q_i}2\right\}.
$$
Hence
$\big(\frac{r_i^+}{q_i}\big)'<\frac1{\theta}\max\{1,\frac2{q_i}\}$
which gives that $w_i^{q_i}\in RH_{\big(\frac{r_i^+}{q_i}\big)'}$. We
have thus shown that
$w^{q_i}\in A_{\frac{q_i}{r_i^-}}\cap
RH_{\big(\frac{r_i^+}{q_i}\big)'}$ which implies that the bilinear
Hilbert transform
is bounded from $L^{q_1}(w_1^{q_1}) \times L^{q_2}(w_2^{q_2})$ into
$L^q(w^q)$.  This completes the proof of \eqref{eqn:bht1-bis}.

Finally, let
$w_i(x)=|x|^{-\frac{a}{q_i}}$ so that
$w(x)=w_1(x)w_2(x)=|x|^{-\frac{a}{q}}$. Then, using the well known
properties of power weights, we have that 
$w_i^{q_i}\in A_{\max\{1,\frac{q_i}2\}}\cap RH_{\max\{1,\frac2{q_i}\}}$ if and only if
$$
1-\max\left\{1,\frac{q_i}2\right\}<a<1
\qquad\text{and}\qquad
-\infty<a<\frac1{\max\{1,\frac2{q_i}\}}
=
\min\left\{1,\frac{q_i}{2}\right\},
$$
and when $\max\{1,\frac{q_i}2\}=1$ we can also allow $a=0$ in the first condition. From all these we easily see that 
\eqref{eq:BH-power} holds provided either $a=0$ or $a$ satisfies \eqref{eq:values-a}. This completes the proof. \qed

\medskip

\subsection*{Proof of Theorem~\ref{thm:bht-vector}} 
The proof of the first part of Theorem ~\ref{thm:bht-vector}
is now straightforward given Theorem~\ref{thm:extrapol-multi}
and Corollary~\ref{cor:extrapol-mult}. Indeed,
Theorem~\ref{thm:bht-DCU-JMM} provides the initial weighted norm
inequalities for the family
\[ \F  = \{ (|B(f,g)|,|f|, |g|) :  f, g \in L^\infty_c. \} \]
 (see the proof of
Theorem~\ref{thm:bht-DCU-JMM}). Thus,
Corollary~\ref{cor:extrapol-mult} applies
and~\eqref{eq:extrap-multi-con-vv} yields~\eqref{eqn:bht-vector1} for
functions $f_k, g_k\in L^\infty_c$. By a standard approximation
argument we get the desired inequality for $f_k\in L^{q_1}(w_1^{q_1})$
and $g_k\in L^{q_2}(w_2^{q_2})$.

\medskip

\medskip 

To prove the second part of Theorem~\ref{thm:bht-vector} we modify the argument
 in the second part of the proof of
Theorem~\ref{thm:bht-DCU-JMM}. Fix $q_i,s_i$ as in the statement;
then by the first part of Theorem \ref{thm:bht-vector} we need to
find $1<p_1,p_2<\infty$ such that
$\frac{1}{p}=\frac{1}{p_1}+\frac{1}{p_2}<1$ and $r_i^-<q_i<r_i^+$ with
\begin{equation}\label{eq:formula-r_ipm-vv}
\frac1{r_i^+}=\frac1{2p_i}
<\frac1{q_i},\frac1{s_i}
<
\frac1{r_i^-}=\frac1{2p_i}+\frac12.
\end{equation}
Since $1<p_1,p_2<\infty$, \eqref{eq:formula-r_ipm-vv} can be rewritten
as 
\begin{equation}\label{eq:formula-r_ipm:equiv-vv}
0\le 2\left(\max\left\{\frac12,\frac1{q_i},\frac1{s_i}\right\}-\frac12\right)
<
\frac1{p_i}
<
2\min\left\{
\frac12,\frac1{q_i},\frac1{s_i}\right\}\le 1.
\end{equation}
Before choosing $p_1,p_2$, we first claim that
\begin{equation}\label{eq:restriction-vv}
\sum_{i=1}^2 \max\left\{\frac12,\frac1{q_i},\frac1{s_i}\right\}<\frac32.
\end{equation}
To show this we argue as we did to prove~\eqref{eq:restriction}: if at least one of the maxima is $\frac12$,
then since the other maxima is strictly smaller than 1 we get
 the desired estimate. If none of the
maxima is $\frac12$, then by the last condition in
\eqref{eq:restriction-exponents-vv},
$$
\sum_{i=1}^2 
\max\left\{\frac12,\frac1{q_i},\frac1{s_i}\right\}
=
\sum_{i=1}^2 \max\left\{\frac1{q_i},\frac1{s_i}\right\} 
<\frac32.
$$

We now choose $p_i$:  fix $\eta_i>0$ and let
\begin{equation}\label{eq:choice-pi-vv}
\frac1{p_i}:=2\left(\max\left\{\frac12,\frac1{q_i},\frac1{s_i}\right\}-\frac12+\eta_i\right),
\qquad i=1,2,
\end{equation}
where we choose the $\eta_i$ sufficiently small so that
\begin{equation}\label{eq:conds-eta-vv:1}
\eta_1+\eta_2<
\frac32-\sum_{i=1}^2 \max\left\{\frac12,\frac1{q_i},\frac1{s_i}\right\}
\end{equation}
and
\begin{equation}\label{eq:conds-eta-vv:2}
0<\eta_i<\min\left\{\frac1{q_i},\frac1{q_i'},\frac1{s_i},\frac1{s_i'},\frac12-\left|\frac1{s_i}-\frac1{q_i}\right|\right\}.
\end{equation}
Such a choice of $\eta_1,\eta_2$ is possible by
\eqref{eq:restriction-vv} and \eqref{eq:restriction-exponents-vv}. By
\eqref{eq:conds-eta-vv:1} we have that
$$
\frac{1}{p}
=
\frac{1}{p_1}+\frac{1}{p_2}
=
-2+2\,\sum_{i=1}^2 \max\left\{\frac12,\frac1{q_i},\frac1{s_i}\right\}+2\sum_{i=1}^2 \eta_i<1.
$$
To prove \eqref{eq:formula-r_ipm:equiv-vv} we first observe that since $\eta_i>0$,
$$
2\left(\max\left\{\frac12,\frac1{q_i},\frac1{s_i}\right\}-\frac12\right)
<
2\left(\max\left\{\frac12,\frac1{q_i},\frac1{s_i}\right\}-\frac12+\eta_i\right)
=
\frac1{p_i}.
$$
To get the second estimate in \eqref{eq:formula-r_ipm:equiv-vv} we
consider two cases. If
$\max\{\frac12,\frac1{q_i},\frac1{s_i}\}=\frac12$, then
$$
\frac1{p_i}
=
2\eta_i
<
2\min\left\{\frac1{q_i},\frac1{s_i}\right\}
=
2\min\left\{
\frac12,
\frac1{q_i},\frac1{s_i}\right\}.
$$
On the other hand, if $\max\{\frac12,\frac1{q_i},\frac1{s_i}\}=\max\{\frac1{q_i},\frac1{s_i}\} $ 
and we write $\frac1{\alpha_i}=\max\{\frac1{q_i},\frac1{s_i}\}$ and $\frac1{\beta_i}=\max\{\frac1{q_i},\frac1{s_i}\}$, we obtain
\begin{align*}
\frac1{p_i}
&=
2\max\left\{\frac1{q_i},\frac1{s_i}\right\}-1+2\eta_i
\\
&<
2\max\left\{\frac1{q_i},\frac1{s_i}\right\}-1+2
\min\left\{\frac1{q_i'},\frac1{s_i'},\frac12-\left|\frac1{s_i}-\frac1{q_i}\right|\right\}
\\
&=
\frac2{\alpha_i}
-1+
2\min\left\{\frac1{\alpha_i'},\frac12-\left(\frac1{\alpha_i}-\frac1{\beta_i}\right)\right\}
\\
&
=
2\min\left\{\frac12,\frac1{\beta_i}\right\}
\\
&
=
2\min\left\{\frac12,\frac1{q_i},\frac1{s_i}\right\}.
\end{align*}
This completes the proof of \eqref{eq:formula-r_ipm:equiv-vv} and hence that of Theorem~\ref{thm:bht-vector}.\qed 

\medskip

\subsection*{Proof of Remark~\ref{remark:iteration}} 
To prove the iterated vector-valued inequality in
Remark~\ref{remark:iteration}, we simply repeat the argument used to
prove the first part of Theorem~\ref{thm:bht-vector}.   For our
starting estimate we form the new
family 
\begin{multline*}
 \F' = \bigg\{ (h,f,g) \\
= \bigg( \bigg(\sum_k |BH(f_k,g_k)|^s\bigg)^{\frac{1}{s}}, 
\bigg(\sum_k |f_k|^{s_1}\bigg)^{\frac{1}{s_1}}, 
\bigg(\sum_k |g_k|^{s_2}\bigg)^{\frac{1}{s_2}} \bigg) : 
f_k, g_k \in L^\infty_c \bigg\}; 
\end{multline*}
then~\eqref{eqn:bht-vector1} gives us the starting estimate
\[ \|h\|_{L^q(w^q)}
\leq C\|f\|_{L^{q_1}(w_1^{q_1})}\|g\|_{L^{q_2}(w_2^{q_2})}. \]
We then  again apply vector-valued extrapolation using the family
\begin{multline*}
 \F''= \bigg\{ (H,F,G) 
= \bigg( \bigg(\sum_j h_j^t\bigg)^{\frac{1}{t}}, 
\bigg(\sum_j  f_j^{t_1}\bigg)^{\frac{1}{t_1}}, 
\bigg(\sum_j g_j^{t_2}\bigg)^{\frac{1}{t_2}} \bigg) : 
(h_j,f_j,g_j) \in \F' \bigg\}
\end{multline*}
to get iterated vector-valued inequalities.   Details are left
to the interested reader. \qed

\medskip

\subsection*{Proof of Corollary~\ref{corol:BHT-vv-power}}

Similar to our approach  in the proof of Corollary~\ref{cor:bht-A1},
here we take a closer
look at the proof of Theorem~\ref{thm:bht-vector}.  Fix
$1<q_1,q_2, s_1, s_2<\infty$ and $w_i^{q_i}$ as in the statement. We
choose $p_i$ as in \eqref{eq:choice-pi-vv}, \eqref{eq:conds-eta-vv:1}
and \eqref{eq:conds-eta-vv:2}, but again  we will choose $\eta_i$ much
smaller.   Then as we proved above, $1<p_1,p_2<\infty$,
$\frac{1}{p}=\frac{1}{p_1}+\frac{1}{p_2}<1$, and
\eqref{eq:formula-r_ipm:equiv-vv} holds. Note that the latter implies
\eqref{eq:formula-r_ipm-vv} and hence, by the first part of
Theorem~\ref{thm:bht-vector}, we obtain that the bilinear Hilbert transform satisfies
\eqref{eqn:bht-vector1}, provided we show that
$w_i^{q_i}\in A_{\frac{q_i}{r_i^-}}\cap
RH_{\big(\frac{r_i^+}{q_i}\big)'}$, where
$$
\frac{q_i}{r_i^-}
=
q_i\left(\frac1{2p_i}+\frac12\right)
=
q_i\left(\max\left\{\frac12,\frac1{q_i},\frac1{s_i}\right\}+\eta_i\right)
=
\max\left\{1,\frac{q_i}2,\frac{q_i}{s_i}\right\}+q_i\eta_i
$$
and
\begin{multline*}
\frac1{\big(\frac{r_i^+}{q_i}\big)'}
=
1-\frac{q_i}{r_i^+}
=
1-\frac{q_i}{2p_i}
=
1-q_i\left(\max\left\{\frac12,\frac1{q_i}, \frac1{s_i},\right\}-\frac12+\eta_i\right)
\\
=
\min\left\{1,\frac{q_i}2,1-q_i\left(\frac1{s_1}-\frac12\right)\right\}-q_i\eta_i.
\end{multline*}
Note that $w_i^{q_i}\in A_{\max\{1,\frac{q_i}2,\frac{q_i}{s_i}\}}$
immediately gives us that $w_i^{q_i}\in A_{\frac{q_i}{r_i^-}}$. On the other
hand, since
$w_i^{q_i}\in RH_{\max\{1,\frac2{q_i},
  [1-q_i(\frac1{s_i}-\frac12)]^{-1}\}}$, by the openness of the
reverse H\"older classes we can find $0<\theta<1$ close to $1$ such that
$w_i^{q_i}\in RH_{\frac1{\theta}\max\{1,\frac2{q_i},
  [1-q_i(\frac1{s_i}-\frac12)]^{-1}\}}$. We therefore assume, in
addition to \eqref{eq:conds-eta-vv:1}, \eqref{eq:conds-eta-vv:2}, that
$0<\eta_i<(1-\theta)\min\big\{\frac12,\frac1{q_i},\frac1{q_i}-\frac1{s_i}+\frac12\big\}$;
this choice is possible because of the first two conditions in
\eqref{eq:restriction-exponents-vv-power-corol}.  But then
\begin{align*}
\frac1{\big(\frac{r_i^+}{q_i}\big)'}
&=
\min\left\{1,\frac{q_i}2,1-q_i\left(\frac1{s_1}-\frac12\right)\right\}-q_i\eta_i
\\
&>
\min\left\{1,\frac{q_i}2,1-q_i\left(\frac1{s_1}-\frac12\right)\right\}
-
q_i(1-\theta)\min\left\{\frac12,\frac1{q_i},\frac1{q_i}-\frac1{s_i}+\frac12\right\}
\\
&=
\theta\min\left\{1,\frac{q_i}2,1-q_i\left(\frac1{s_i}-\frac12\right)\right\}.
\end{align*}
Hence
$\big(\frac{r_i^+}{q_i}\big)'<\frac1{\theta}\max\{1,\frac2{q_i},
[1-q_i(\frac1{s_i}-\frac12)]^{-1}\}$, so
$w_i^{q_i}\in RH_{\big(\frac{r_i^+}{q_i}\big)'}$. We have thus shown
that
$w^{q_i}\in A_{\frac{q_i}{r_i^-}}\cap
RH_{\big(\frac{r_i^+}{q_i}\big)'}$ which  yields
\eqref{eqn:bht-vector1}.

To complete our proof we need to establish
\eqref{eqn:bht-vector1:pwerw-wts}.  Let
$w_i(x)=|x|^{-\frac{a}{q_i}}$ so that
$w(x)=w_1(x)w_2(x)=|x|^{-\frac{a}{q}}$. Then, using the well known
properties of power weights, we have that 
$w_i^{q_i}\in A_{\max\{1,\frac{q_i}2,\frac{q_i}{s_i}\}}\cap RH_{\max\{1,\frac2{q_i}, [1-q_i(\frac1{s_i}-\frac12)]^{-1}\}}$ if and only if
$$
1-\max\left\{1,\frac{q_i}2,\frac{q_i}{s_i}\right\}<a<1;
$$
when $\max\{1,\frac{q_i}2\}=1$ we can also allow $a=0$, and
$$
-\infty<a<\frac1{\max\{1,\frac2{q_i}, [1-q_i(\frac1{s_i}-\frac12)]^{-1}\}}
=
\min\left\{1,\frac{q_i}2, 1-q_i\left(\frac1{s_i}-\frac12\right)\right\}
.
$$
From all these estimates we  see that
\eqref{eqn:bht-vector1:pwerw-wts} holds provided 
$a\in\{0\}\cup(a_-,a_+)$ with $a_\pm$ defined in
\eqref{eq:rest-a-v-v}. This completes the proof. \qed

\subsection*{Proof of Theorem~\ref{thm:CZ-MZ-new}}
The desired result follows directly from extrapolation.  
Fix $1<r<2$ and define the family of $(m+1)$-tuples
\begin{multline*}
\F = \bigg\{ (F,F_1,\ldots,F_m) \\
= \bigg( \bigg(\sum_{k_1,\ldots,k_m}
|T(f^1_{k_1},\ldots,f^m_{k_m})|^r\bigg)^{\frac{1}{r}},
\bigg(\sum_{k_1} |f^1_{k_1}|^r\bigg)^{\frac{1}{r}},\ldots,
\bigg(\sum_{k_m} |f^m_{k_m}|^r\bigg)^{\frac{1}{r}}\bigg)
: f_{k_j}^j \in L^\infty_c \bigg\}.\!\!\!\null 
\end{multline*}
Now fix $1<q_1,\ldots,q_m<r<2$ and let $\frac{1}{q}=\sum \frac{1}{q_j}$.
Then by Theorem~\ref{theorem:CZ-MZ}, for all weights $w_j$ such that
weights $w_j^{q_i} \in A_{q_j}$, and $(F,F_1,\ldots,F_m) \in
\F$,
\begin{equation} \label{eqn:CM-MZ-base}
 \|F\|_{L^q(w^q)} \leq
C\prod_{j=1}^m \|F_j\|_{L^{q_j}(w_j^{q_j})}.  
\end{equation}
Therefore, by Corollary~\ref{cor:extrapol-mult} applied with $r_j^- =
1$, $r_j^+=\infty$, $1\leq j \leq m$, we immediately conclude that for
any $1<q_1,\ldots,q_m<\infty$ and weights $w_j^{q_i} \in A_{q_j}$,
inequality \eqref{eqn:CM-MZ-base} holds, which 
yields~\eqref{eqn:CZ-MZ1} for functions in $L^\infty_c$.  The desired
inequality then follows for $f_{k_j}^j \in L^{q_j}(w_j^{q_j})$ by a standard
approximation argument. \qed

\medskip

\section{More general vector-valued inequalities}
\label{section:generalize}

In this section we explain how to obtain, via extrapolation,
vector-valued inequalities in a larger range than we proved
in Theorem~\ref{thm:bht-vector}.   The starting point is  implicit
in the proof of \cite[Corollary 4]{Culiuc:2016wr}:  from it one can show that
\eqref{eqn:bht0} holds provided
\begin{equation}\label{eq:weights-new}
w_i^{p_i}\in A_{1+(1-\theta_i)(p_i-1)}\cap RH_{\frac1{1-\theta_3}},      
\end{equation}
where $1<p_1,\,p_2,\,p<\infty$ with $\frac1p_1+\frac1p_2=\frac1p$, and where 
$\theta_1,\theta_2,\theta_3\in(0,1)$ are arbitrary parameters satisfying
\begin{equation}\label{p1-2-3:conds}
\frac{\theta_1}{p_1'}\le\frac12,
\qquad
\frac{\theta_2}{p_2'}\le \frac12,
\qquad
\frac{\theta_3}{p}\le\frac12,
\qquad
\frac{\theta_1}{p_1'}+\frac{\theta_2}{p_2'}+\frac{\theta_3}{p}=1.
\end{equation}
In \cite{Culiuc:2016wr} the authors chose
$\theta_1=\theta_2=\theta_3=\frac12$, which then gives Theorem \ref{thm:bht-cdo}.

If we now fix the parameters $\theta_1,\theta_2,\theta_3\in(0,1)$, we can rewrite \eqref{eq:weights-new} as 
\begin{equation}\label{eq:new-values-ri-pm}
w_i^{p_i}\in A_{\frac{p_i}{r_i^-}}\cap RH_{\big(\frac{r_i^+}{p_i}\big)'},
\qquad\text{where}\quad
\frac1{r_i^-}=1-\frac{\theta_i}{p_i'}
\quad\text{and}\quad
\frac1{r_i^+}=\frac{\theta_3}{p_i}.
\end{equation}
Given this, we can apply our extrapolation result to obtain
vector-valued inequalities by varying $p_1,p_2,p$ and
$\theta_1,\theta_2,\theta_3$. We claim that, as a result,
\eqref{eqn:bht-vector1} holds (taking $w_1=w_2\equiv 1$ for
simplicity, but of course some natural weighted norm inequalities are
also possible) whenever $1<s_1,s_2,q_1,q_2<\infty$,
$\frac{1}{s}=\frac{1}{s_1}+\frac{1}{s_2}<\frac32$,
$\frac{1}{q}=\frac{1}{q_1}+\frac{1}{q_2}<\frac{3}{2}$ and if there
exist $0\le \gamma_1,\gamma_2,\gamma_3<1$ with
$\gamma_1+\gamma_2+\gamma_3=1$ such that
\begin{equation}\label{eq:new-cond-vv:1}
\max\left\{\frac{1}{s_1},\frac{1}{q_1}\right\}<\frac{1+\gamma_1}{2},
\ \
\max\left\{\frac{1}{s_2},\frac{1}{q_2}\right\}<\frac{1+\gamma_2}{2},
\ \ 
\max\left\{\frac{1}{s'},\frac{1}{p'}\right\}<\frac{1+\gamma_3}{2},
\end{equation}
and, additionally,
\begin{equation}\label{eq:new-cond-vv:2}
\min\left\{\frac{1}{s_1},\frac{1}{q_1}\right\} 
+
\min\left\{\frac{1}{s_2},\frac{1}{q_2}\right\}
>
\frac{1-\gamma_3}{2}.
\end{equation}
If we compare our conditions with those in
\cite[Theorem~5]{Benea:2016th} (see also
\cite[Appendix~A]{Culiuc:2016wr}), we see that ours impose the extra
restrictions \eqref{eq:new-cond-vv:2} and $s_i, q_i<\infty$.  Also,
note that the last condition in \eqref{eq:new-cond-vv:1} is implied by
\eqref{eq:new-cond-vv:2}; nevertheless we make it explicit in order to
compare our conditions with those of \cite[Theorem~5]{Benea:2016th}.

We now sketch how to prove our claim.   Define
$$
m_1=\min\left\{\frac{1}{s_1},\frac{1}{q_1}\right\},
\quad
m_2=\min\left\{\frac{1}{s_2},\frac{1}{q_2}\right\},
\quad
\widetilde{m}_1= \frac2{1-\gamma_3}m_1,
\quad
\widetilde{m}_2= \frac2{1-\gamma_3}m_2.
$$
With this notation, \eqref{eq:new-cond-vv:2} becomes
$\widetilde{m}_1+\widetilde{m}_2>1$.  The first step is to show that
there exist $0<\eta_1,\eta_2<1$ such that
\begin{equation}\label{eq:choice-etai}
\eta_1+\eta_2=1,
\qquad
\eta_1<\widetilde{m}_1,
\qquad
\eta_2<\widetilde{m}_2.
\end{equation}
To prove this we consider two cases. If
$|\widetilde{m}_1-\widetilde{m}_2|<1$, we just need to pick
$\eta_1:=\frac12+\frac{\widetilde{m}_1-\widetilde{m}_2}{2}$,
$\eta_2:=\frac12+\frac{\widetilde{m}_2-\widetilde{m}_1}{2}$. On the
other hand, if $|\widetilde{m}_1-\widetilde{m}_2|\ge 1$ then either
$\widetilde{m}_1\ge 1$ or $\widetilde{m}_2\ge 1$. If
$\widetilde{m}_1\ge 1$,  let $\eta_1=1-\epsilon$, $\eta_2=\epsilon$ with
$0<\epsilon\ll 1$;  if $\widetilde{m}_2\ge 1$, let $\eta_1=\epsilon$,
$\eta_2=1-\epsilon$ with $0<\epsilon\ll 1$. 

Once $\eta_1,\eta_2$ are
chosen we consider two cases. When $0<\eta_1\le\eta_2<1$, we take
$$
\frac{2\eta_2}{\eta_1}
<
p_1
<
\frac{2}{(1-\gamma_3)\eta_1},
\qquad
p_2=p_1\frac{\eta_1}{\eta_2},
\qquad
p=p_1\eta_1=p_2\eta_2.
$$
Then we have that $p_1,p_2>2$ and 
\begin{equation}\label{eq:relat-p:1}
1= \eta_1+\eta_2 \le 2\eta_2<p_1\eta_1=p<\frac{2}{1-\gamma_3}.
\end{equation}

When $0<\eta_2<\eta_1<1$, we choose
$$
\frac{2\eta_1}{\eta_2}
<
p_2
<
\frac{2}{(1-\gamma_3)\eta_2},
\qquad
p_1=p_2\frac{\eta_2}{\eta_1},
\qquad
p=p_1\eta_1=p_2\eta_2.
$$
Again we have $p_1,p_2>2$ and 
\begin{equation}\label{eq:relat-p:2}
1= \eta_1+\eta_2 < 2\eta_1<p_2\eta_2=p<\frac{2}{1-\gamma_3}.
\end{equation}

In both cases we have
$
\frac1{p_1}+\frac1{p_2}=\frac1p(\eta_1+\eta_2)=\frac1p<1.
$
Now let
$$
\theta_1=p_1'	\frac{1-\gamma_1}{2},
\qquad
\theta_2=p_2'\frac{1-\gamma_2}{2},
\qquad
\theta_3=p\frac{1-\gamma_3}{2}.
$$
Then $\theta_1,\theta_2,\theta_3>0$ since
$\gamma_1,\gamma_2,\gamma_3<1$. From \eqref{eq:relat-p:1} or
\eqref{eq:relat-p:2} we have that $\theta_3<1$, and, since
$p_1,p_2>2$, it follows that $\theta_i<1-\gamma_i\le 1$ for
$i=1,2$. We also have that $\frac{\theta_1}{p_1'}\le\frac12$,
$\frac{\theta_2}{p_2'}\le \frac12$, and $\frac{\theta_3}{p}\le\frac12$
since $\gamma_1,\gamma_2,\gamma_3 \ge 0$. Finally, since we assumed
that $\gamma_1+\gamma_2+\gamma_3=1$, we get
$$
\frac{\theta_1}{p_1'}+\frac{\theta_2}{p_2'}+\frac{\theta_3}{p}
=
\frac{1-\gamma_1}{2}+\frac{1-\gamma_2}{2}+\frac{1-\gamma_3}{2}
=1.
$$
Therefore, \eqref{p1-2-3:conds} holds and, as observed above, 
it follows that \eqref{eq:weights-new} yields \eqref{eqn:bht0}.

By extrapolation (arguing as we did in the proof of Theorem
\ref{thm:bht-vector}) and using \eqref{eq:new-values-ri-pm},  we get that \eqref{eqn:bht-vector1} holds
provided 
\begin{equation}\label{eq:needed-finish}
\frac{\theta_3}{p_i}
=
\frac1{r_i^+}.
<\frac1{s_i},\frac1{q_i}<
\frac1{r_i^-}
=
1-\frac{\theta_i}{p_i'},
\qquad
i=1,2.
\end{equation}
Hence, we need to show that \eqref{eq:new-cond-vv:1} and
\eqref{eq:new-cond-vv:2} imply that~\eqref{eq:needed-finish} holds.  First,
from~\eqref{eq:new-cond-vv:1} we get that
$$
\max\left\{\frac{1}{s_i},\frac{1}{q_i}\right\}<\frac{1+\gamma_i}{2}
=
1-\frac{\theta_i}{p_i'}
=
\frac1{r_i^-}
,\qquad i=1,2.
$$
Second, \eqref{eq:choice-etai} yields
$$
\frac1{r_i^+}
=
\frac{\theta_3}{p_i}
=
\frac{p}{p_i}\frac{1-\gamma_3}{2}
=
\eta_i
<
\widetilde{m}_i\frac{1-\gamma_3}{2}
=
m_i
=
\min\left\{\frac{1}{s_i},\frac{1}{q_i}\right\}.
$$
Hence, \eqref{eq:needed-finish} holds, and this completes our sketch
of the proof.

\bibliographystyle{plain}
\bibliography{off-diag-extrapol}

\end{document}